\documentclass[finall]{autart}    
\newcommand{\EQ}{\begin{eqnarray}}
\newcommand{\EN}{\end{eqnarray}}
\newcommand{\EQQ}{\begin{eqnarray*}}
\newcommand{\ENN}{\end{eqnarray*}}
\newcommand{\nnum}{\nonumber}
\renewcommand{\theequation}{\thesection.\arabic{equation}}

\newcommand{\bremark}{\begin{remark} \begin{rm} }
\newcommand{\eremark}{ \end{rm}  \qed\\
\end{remark} }
\newcommand{\btheorem}{\begin{thm} \begin{rm} }
\newcommand{\etheorem}{ \end{rm}  \qed\\
\end{thm} }
\newcommand{\blemma}{\begin{lemma} \begin{rm} }
\newcommand{\elemma}{ \end{rm}  \qed\\
\end{lemma} }
\newcommand{\bcorollary}{\begin{corollary} \begin{rm} }
\newcommand{\ecorollary}{ \end{rm}  \qed\\
\end{corollary} }
\usepackage{graphicx}          
   \usepackage{cite}                        
\usepackage{amssymb}                              
\usepackage{dsfont}
\usepackage{eufrak}  
\usepackage{mathabx} 
\begin{document}

\begin{frontmatter}

\title{ Averaging for nonlinear systems evolving on Riemannian manifolds}                                                

\thanks[footnoteinfo]{This research was partially supported by the Australian Research Council.}

\author[First]{Farzin Taringoo}\ead{ftaringoo@unimelb.edu.au},   
\author[First]{Dragan Ne\v{s}i\'{c}}\ead{d.nesic@unimelb.edu.au}, 
\author[First]{Ying Tan}\ead{yingt@unimelb.edu.au},
\author[First]{Peter M. Dower}\ead{pdower@unimelb.edu.au}              

\address[First]{Electrical and Electronic Engineering Department, The University of
Melbourne, Victoria, 3010, Australia}  


\begin{abstract}                          
This  paper presents an averaging method  for  nonlinear systems defined on Riemannian manifolds. We extend closeness of solutions  results  for ordinary differential equations on $\mathds{R}^{n}$ to dynamical systems defined on Riemannian manifolds by employing differential geometry. A generalization of closeness of solutions for periodic dynamical systems on compact time intervals is derived for dynamical systems evolving on compact Riemannian manifolds. Under local asymptotic (exponential) stability of the average vector field, we further relax the compactness of the ambient Riemannian manifold and obtain the closeness of solutions on the infinite time interval by employing the notion of uniform normal neighborhoods of an equilibrium point of a vector field. These results are also presented for time-varying dynamical systems where their averaged systems are almost globally asymptotically or exponentially stable on compact manifolds. The main results of the paper are illustrated by several examples.\end{abstract}

\end{frontmatter}

\section{Introduction}

  Perturbation theory is a class of mathematical methods used to find approximations of solutions of dynamical systems which cannot be solved directly, see \cite{Kha, sand, Holmes}. Averaging is a powerful perturbation based tool that has applications in the study of time-varying linear and nonlinear dynamical systems. Where applicable, averaging can provide closeness of solutions results for the solutions of trajectories of such systems related to those of a corresponding averaged system. As the trajectories of this averaged system can be substantially simpler than those of the original time-varying system, stability analysis can be simplified by exploiting closeness of solutions results provided by the averaging, see \cite{nes3, nes4,Bullo,volo, Bai, Per}. Averaging results have been developed for numerous classes of dynamical systems and differential inclusions (see  \cite{Bit,Don, nes1, nes2, Sus}) including dynamical systems on Lie groups, see \cite{Leo1, Leo2, Leo3, Gur}.

The state spaces of many dynamical systems constitute  Riemannian manifolds (see \cite{Lewis, Bloch, Arnold, mar1, Sastry}) and consequently their analyses require differential geometric tools. Examples of such systems can be found in many mechanical settings, see \cite{Lewis, Bloch}. 
In this paper, averaging is extended to a particular class of dynamical systems evolving on Riemannian manifolds.   Such systems arise naturally in classical mechanics (see \cite{Lewis, Bloch, Arnold}) where the state space of the dynamical system is restricted to such a manifold. A version of averaging methods for dynamical systems on Lie groups is introduced in \cite{Leo1, Leo2, Leo3}. We address the problem of closeness of solutions on  finite and infinite time horizons on Riemannian manifolds. These results generalize those  presented in \cite{Kha}, Chapter 10. In the case of compact time intervals, the analyses are presented for dynamical systems on compact Riemannian manifolds. 

By employing the notion of \textit{Levi-Civita} connection on Riemannian manifolds, we study the closeness of solutions of vector fields where the closeness is exploited with respect to the Riemannian distance function, see \cite{Lee4}. Using a version of stability theory for systems evolving on Riemannian manifolds (see \cite{forni,Ang, Lewis}) we extend the closeness of solutions results (\cite{Kha, nes, sand}) to the infinite time interval where average systems are assumed to be locally asymptotically or exponentially stable. We use the scaling technique to bound the Riemannian metric by the Euclidean one (see \cite{Lee3, jost, Pet}) on a precompact (see \cite{Lee3}) neighborhood of an equilibrium  of the average system  in its uniform normal neighborhood and invoke some of the standard results of the stability theory presented in \cite{Kha}. Geometric features of the normal neighborhoods such as existence of unique length minimizing geodesics and their local representations enable us to closely relate the results obtained for dynamical systems in $\mathds{R}^{n}$ to those in Riemannian manifolds.

In terms of exposition,  Section 2  presents some  mathematical preliminaries needed for the analyses of the paper. Section 3 presents the main averaging results for dynamical systems on Riemannian manifolds on finite time horizon together with some numerical examples.  The results of Section 3 are strengthened to the infinite time horizon limit in Sections 4 and 5 by employing a notion of stability on Riemannian manifolds. 
\section{Preliminaries}In this section we provide the differential geometric material which is necessary for the analyses presented in the rest of the paper. Table I summarizes key notation used throughout:\\
\begin{table}[ht]
\caption{Notation and descriptions} 
\centering 
\begin{tabular}{c c} 
\hline\hline 
Symbol& Description  \\ [0.5ex] 
\hline 
$M$ & Riemannian manifold \\ 
$\mathfrak{X}(M)$ & space of smooth time-invariant\\ &vector fields on $M$\\ 
$\mathfrak{X}(M\times \mathds{R})$ & space of smooth time-varying\\ &vector fields on $M$\\ 
$\mathfrak{X}(\mathds{R}\times M)$ & space of smooth parameter-varying\\ &vector fields on $M$\\
$C^{\infty}(M)$ & space of smooth functions on $M$\\
$T_{x}M$ & tangent space at $x\in M$  \\
$T^{*}_{x}M$ & cotangent space at $x\in M$  \\
$TM$ & tangent bundle of $M$\\
$T^{*}M$ & cotangent bundle of $M$\\
$\frac{\partial}{\partial x_{i}}$ & basis tangent vectors at $x\in M$\\
$dx_{i}$ & basis cotangent vectors at $x\in M$\\
$f(x,t)$ & time-varying vector fields on $M$  \\
$||f||_{g}$ & Riemannian norm of $f$\\
$||f||_{e}$ & Euclidean norm of $f$\\
$g(.,.)$ & Riemannian metric on $M$ \\ 
$d(.,.)$ & Riemannian distance on $M$\\
$\nabla$ & (Levi-Civita) Connection on $M$\\
$\Phi_{f}$ & flow associated with $f$\\
$T\Phi_{f}$ & push-forward of $\Phi_{f}$\\
$T^{*}\Phi_{f}$ & pull-back of $\Phi_{f}$\\
$\mathds{R}_{>0}$& $(0,\infty)$\\
$\mathds{R}_{\geq 0}$& $[0,\infty)$\\

\hline 
\end{tabular}
\label{table:nonlin} 
\end{table}
\newtheorem{definition}{Definition}
\begin{definition}
A Riemannian manifold $M$ is a differentiable manifold together with a Riemannian metric $g, \hspace{.2cm}x\in M$, where $g:T_{x}M\times T_{x}M\rightarrow \mathds{R}$ is symmetric and positive definite where $T_{x}M$ is the tangent space at $x\in M$ (see \cite{Lee2}, Chapter 3). For $M=\mathds{R}^{n}$, the Riemannian metric $g$ is given by 
 \EQ g\left(\frac{\partial}{\partial x_{i}},\frac{\partial}{\partial x_{j}}\right)=\delta_{ij},\hspace{.2cm}i,j=1,...,n,\EN
 where $\delta_{ij}$ is the Kronecker delta. 
\end{definition}
\begin{definition}
For a given smooth mapping $F:M\rightarrow N$ from manifold $M$ to manifold $N$ the pushforward $TF$ is defined as a generalization of the Jacobian of smooth maps in Euclidean spaces as follows:
\EQ TF:T_{x}M\rightarrow T_{F(x)}N,\EN
where
\EQ TF(X_{x})\circ f=X_{x}(f\circ F),\hspace{.2cm}X_{x}\in T_{x}M, f\in C^{\infty}(N).\nnum\\\EN
The pullback $T^{*}F$ is defined by
\EQ T^{*}F:T^{*}_{F(x)}N\rightarrow T^{*}_{x}M,\EN
where
\EQ T^{*}F(\omega)\circ X_{x}=\omega(TF(X_{x})),\hspace{.2cm}X_{x}\in T_{x}M, \omega\in T^{*}_{F(x)}N,\nnum\\\EN
where $T^{*}M$ is the cotangent bundle of $M$ (see \cite{Lee2} Chapters 3 and 6).
\end{definition}
In this paper we restrict the analysis to connected finite dimensional Riemannian manifolds.
On an $n$ dimensional Riemannian manifold $M$, the length function of a smooth curve $\gamma:[a,b]\rightarrow M$ is defined as follows:
\EQ \ell(\gamma)=\int^{b}_{a}\big(g(\dot{\gamma}(t),\dot{\gamma}(t))\big)^{\frac{1}{2}}dt,\EN
in which $g$ denotes the Riemannian metric on $M$.
Consequently  we can define a metric (distance) $d$ on an $n$ dimensional Riemannian manifold $(M,g)$ as follows:
\EQ\label{l} &&d:M\times M\rightarrow \mathds{R},\hspace{.2cm}\nnum\\&& d(x,y)=\inf_{\gamma:[a,b]\rightarrow M}\int^{b}_{a}\big(g(\dot{\gamma}(t),\dot{\gamma}(t))\big)^{\frac{1}{2}}dt,\nnum\\ \EN where $\gamma:[a,b]\rightarrow M$ is a piecewise smooth path and $\gamma(a)=x,\gamma(b)=y$.
For $M=\mathds{R}^{n}$ we have
\EQ d(x,y)=\left(\sum^{n}_{i=1}(x_{i}-y_{i})^{2}\right)^{\frac{1}{2}}.\EN

The following theorem ensures that for any connected Riemannian manifold $M$, any pair of points $x,y\in M$ can be connected by a piecewise smooth path $\gamma$. This notion is used to construct a family of curves in the proof of one of the main results of the paper.\\
\newtheorem{theorem}{Theorem}
\begin{theorem}[\hspace{-.025cm}\cite{Lee3}, Page 94]
\label{t1}
Suppose $(M,g)$ is an $n$ dimensional connected Riemannian manifold. Then, for any pair $p,q\in M$, there exists a piecewise smooth path which connects $p$ to $q$.
\end{theorem}
Employing the distance function above it can be shown that $(M,d)$ is a metric space. This is formalized by the next theorem.\\
\begin{theorem}[\hspace{-.02cm}\cite{Lee3}, Page 94]
\label{t2}
With the distance function $d$ defined in (\ref{l}), any connected Riemannian manifold is a metric space where the induced topology is same as the manifold topology. 
\end{theorem}
For a smooth $n$ dimensional Riemannian manifold $M$, a \textit{linear connection} is defined by the following map (see \cite{Lee3})
\EQ \nabla:TM\times TM\rightarrow TM,\EN
where for all $a,b\in \mathds{R}$, $f,h\in C^{\infty}(M)$ and $X,Y,Z\in \mathfrak{X}(M)$, 
\EQ \label{kos1}\nabla_{fX+hY}Z=f\nabla_{X}Z+h\nabla_{Y}Z,\EN
\EQ \nabla_{X}(aY+bZ)=a\nabla_{X}Y+b\nabla_{X}Z,\EN
\EQ \label{lev}\nabla_{X}(fY)=f\nabla_{X}Y+X(f)Y.\EN
The \textit{Levi-Civita} connection $\nabla:\mathfrak{X}(M)\times \mathfrak{X}(M)\rightarrow\mathfrak{X}(M)$ is the unique linear connection on $M$ (see \cite{Lee3}, Theorem 5.4)  which is torsion free and compatible with the Riemannian metric $g$ as follows:  
\EQ \label{kir3}&&\mbox{compatibility with}\hspace{.2cm} g\nnum\\&&Xg(Y,Z)=g(\nabla_{X}Y,Z)+g(Y,\nabla_{X}Z),\EN
\EQ  \label{levi} \hspace{-.5cm}&&(i)(\mbox{torsion free}): \nabla_{X}Y-\nabla_{Y}X=[X,Y],\hspace{.5cm}\nnum\\\hspace{-.2cm}&&(ii): \nabla_{X}f=X(f),\EN
where 
\EQ [X,Y](f)=X(Y(f))-Y(X(f)).\EN

\begin{figure}
\begin{center}
\vspace*{-1.25cm}\hspace*{-.25cm}\includegraphics[scale=.28]{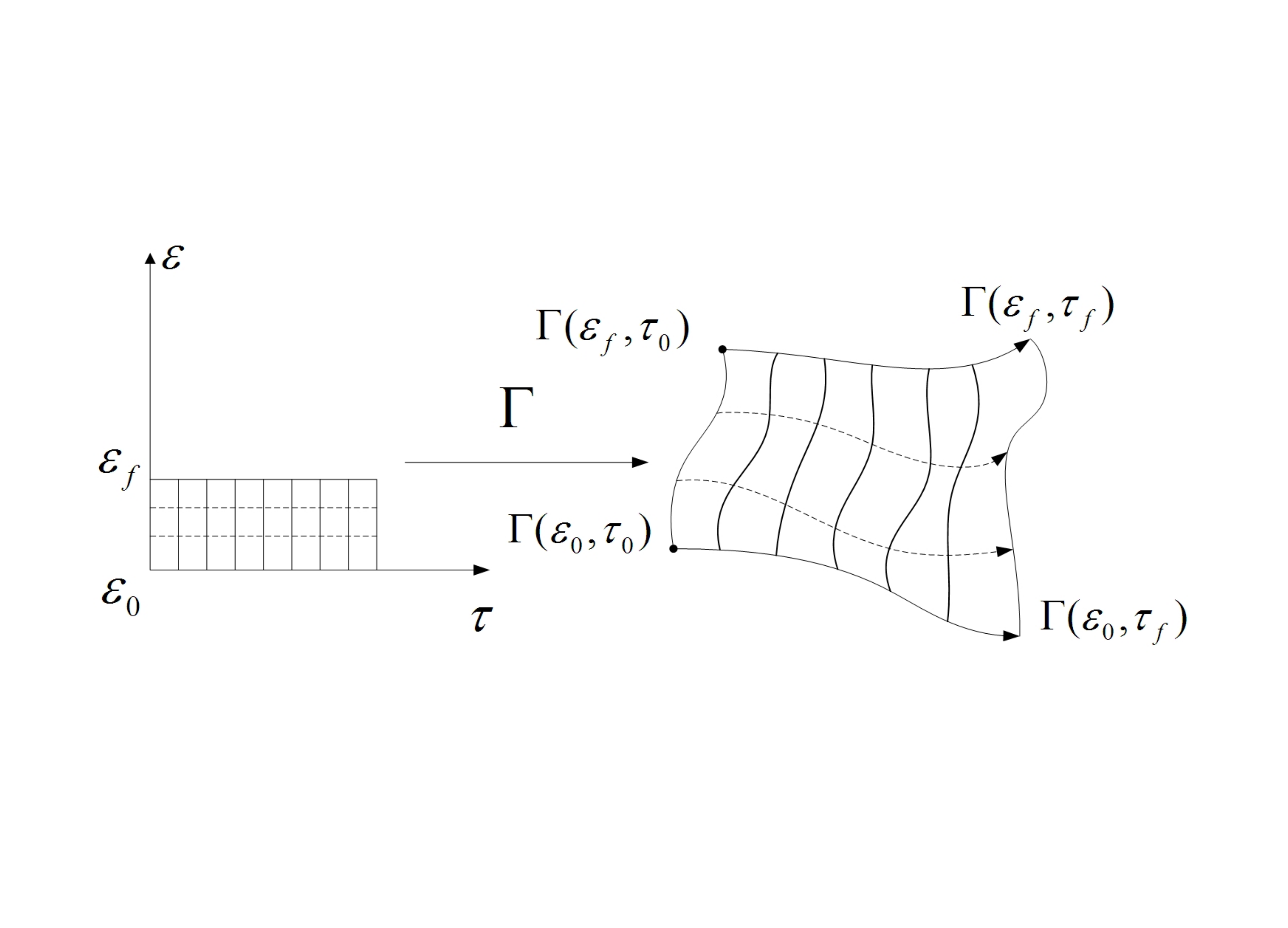}
 \vspace*{-2.25cm}\caption{Admissible family of curves (Definition \ref{kir})}
     \label{ff2}
      \end{center}
   \end{figure}

For $M=\mathds{R}^{n}$ we have
\EQ \nabla_{\frac{\partial}{\partial x_{i}}}\frac{\partial}{\partial x_{j}}=0,\hspace{.2cm}i,j=1,...,n.\EN
\begin{definition}[\hspace{-.02cm}\cite{Lee3}, Page 96]
\label{kir}
An \textit{admissible family} of curves on $M$ is a continuous map \\$\Gamma:(\epsilon_{0},\epsilon_{f})\times [\tau_{0},\tau_{f}]\rightarrow M,\hspace{.2cm}\epsilon\in(\epsilon_{0},\epsilon_{f}) ,\tau\in [\tau_{0},\tau_{f}],\\ \epsilon_{0},\epsilon_{f},\tau_{0},\tau_{f}\in \mathds{R}$ such that $\Gamma$ is smooth with respect to $\epsilon$ and $\tau$ (see Figure \ref{ff2}). 
\end{definition}

Let us denote the tangent vectors obtained by differentiating $\Gamma$ with respect to $\epsilon$ and $\tau$ by
\EQ \partial_{\tau}\Gamma(\epsilon,\tau)\doteq\frac{\partial}{\partial\tau}\Gamma(\epsilon,\tau),\hspace{.2cm}\partial_{\epsilon}\Gamma(\epsilon,\tau)\doteq\frac{\partial}{\partial\epsilon}\Gamma(\epsilon,\tau).\EN
Note that in general $\partial_{\tau}\Gamma(\epsilon,\tau)$ and $\partial_{\epsilon}\Gamma(\epsilon,\tau)$ do not necessarily define vector fields on $M$ since the image of $\Gamma$ may not cover $M$. However, the following lemma enables us to employ the Levi-Civita connection $\nabla$ of $M$ in order to analyze the variation of  $\partial_{\tau}\Gamma(\epsilon,\tau)$ and $\partial_{\epsilon}\Gamma(\epsilon,\tau)$ with respect to vector fields on $M$. 
\newtheorem{lemma}{Lemma}
\begin{lemma}[\hspace{-.02cm}\cite{Lee3}, Page 50, Lemma 4.1]
\label{l1}
Consider $\gamma:(-\epsilon,\epsilon)\rightarrow M,\hspace{.2cm}\epsilon\in \mathds{R}_{>0}$ such that $\gamma(0)=p\in M$ and $\dot{\gamma}(0)=X_{p}\in T_{p}M$. If two vector fields $Y$ and $\tilde{Y}$ agree along $\gamma$, then 
\EQ\nabla_{X_{p}}Y|_{p}= \nabla_{X_{p}}\tilde{Y}|_{p}.\EN\hspace*{7.8cm}\qed
\end{lemma}

Note that by (\ref{levi}), $\nabla$ is torsion free, i.e. $\nabla_{X}Y-\nabla_{Y}X=[X,Y]$. Also note that $[\frac{\partial}{\partial \epsilon},\frac{\partial}{\partial \tau}]=0$, then we have
\EQ \label{kir4}\nabla_{\frac{\partial}{\partial \epsilon}}\partial_{\tau}\Gamma(\epsilon,\tau)=\nabla_{\frac{\partial}{\partial \tau}}\partial_{\epsilon}\Gamma(\epsilon,\tau).\EN
The property above will be used to extend standard averaging techniques to dynamical systems defined on Riemannian manifolds.
In particular, this paper focuses on dynamical systems governed by differential equations on $M$ defined by
\EQ &&\hspace{-0cm}\dot{x}(t)=f(x(t),t),\hspace{.2cm} \nnum\\&&\hspace{-0cm}f(x(t),t)\in T_{x(t)}M,\hspace{.2cm} x(0)=x_{0}\in M, t\in[t_{0},t_{f}],\EN
where $x(t)$ denotes the state at time $t\in[t_0, t_f]$.
The time dependent flow associated with a differentiable time dependent vector field $f$ is a map $\Phi_{f(x,t)}$ satisfying 
\EQ  &&\Phi_{f}:[t_0, t_{f}]\times [t_{0}, t_{f}]\times M\rightarrow M, \nnum\\&& (t_{0},s,x)\mapsto \Phi_{f}(s,t_{0},x)\in M,\EN
and
\EQ \left.\frac{d\Phi_{f}(s,t_{0},x)}{ds}\right|_{s=t}=f(x(t),t)\in T_{x(t)}M.\EN 
One may show that, for a smooth vector field $f$, the integral flow $\Phi_{f}(s,t_{0},.)$ is a local diffeomorphism, see \cite{Lee2}.
  In this paper, on non compact manifolds, we assume that the vector field $f$ is smooth and \textit{complete}, i.e. $\Phi_{f}$ exists for all $t\in [t_{0},\infty)$.
  \subsection{Geodesic Curves}
   As known  geodesics are defined as length minimizing curves on Riemannian manifolds \cite{jost}. The solution of the Euler-Lagrange variational problem associated with the length minimizing problem shows that  all geodesics on $M$  must locally satisfy  the system of ordinary differential equations given by (see \cite{Lee3}, Theorem 4.10)
\EQ \label{geo}\ddot{x}_{i}(t)+\sum^{n}_{j,k=1}\Gamma^{i}_{j,k}\dot{x}_{j}(t)\dot{x}_{k}(t)=0,\quad i=1,...,n,\EN
where
\EQ\label{cris} \Gamma^{i}_{j,k}\doteq\frac{1}{2}\sum^{n}_{l=1}g^{il}(g_{jl,k}+g_{kl,j}-g_{jk,l}),\quad g_{jl,k}\doteq\frac{\partial}{\partial x_{k}}g_{jl},\nnum\\\EN
in which $g$ denotes the Riemannian metric on $M$, and $i,j,k\in[1,...,n]$, $n\doteq \dim(M)$. Note
that $[g^{ij}] = [g_{ij}]^{-1}$.
 \begin{definition}[\hspace{-.02cm}\cite{Lee3}, Page 72]
 The  restricted exponential map is defined by 
 $ \exp_{x}:T_{x}M\rightarrow M,\hspace{.2cm}\exp_{x}(v)=\gamma_{v}(1), v\in T_{x}M,$
 where $\gamma_v:[0,1]\rightarrow M$ is the unique maximal geodesic satisfying $\gamma_v(0)=x\in M$,
$\dot\gamma_v(0) = v\in T_x M$, see \cite{Lee3}, Theorem 4.10.
 \end{definition}
 For the economy of notation, in this paper we refer the restricted exponential maps as exponential maps. 
 For $x\in M$, consider a $\delta$ ball in $T_{x}M$ such that $B_{\delta}(0)\doteq\{v\in T_{x}M \hspace{.1cm}|\hspace{.1cm} ||v||_{g}<\delta\}$. Then the geodesic ball is defined by the following definition.
\begin{lemma}[\hspace{-.02cm}\cite{Lee3}]
 \label{eun}
 For any $x\in M$ there exists a neighborhood $B_{\delta}(0)$ in $T_{x}M$ on which $\exp_{x}$ is a diffeomorphism. 
 \end{lemma}
 \begin{definition}[\hspace{-.02cm}\cite{Lee3}]
 In a neighbourhood of  $x\in M$ where $\exp$ is a local diffeomorphism (this neighborhood always exits by Lemma \ref{eun}), a geodesic ball of radius $0<\delta$ is $\exp_{x}(B_{\delta}(0))\subset M$. Also, we call $\exp_{x}(\overline{B}_{\delta}(0))$ a closed geodesic ball of radius $\delta$. 
 \end{definition}
 \begin{definition}
 For a vector space $V$, a \textit{star-shaped neighborhood} of $0\in V$ is any open set $U$ such that if $u\in U$ then $\alpha u\in U, \alpha\in[0,1]$.
 \end{definition}
 \begin{definition}[\hspace{-.02cm}\cite{Lee3}]
 A normal neighborhood around $x\in M$ is any open neighborhood of $x$ which is a diffeomorphic preimage of a star shaped neighborhood of $0\in T_{x}M$ under $\exp$ map. A uniform normal neighborhood of $x$ is any open set which is contained in a geodesic ball of radius $\delta>0$ for all its points.  
 \end{definition}

 \begin{lemma}[\hspace{-.02cm}\cite{Lee3}]
 \label{un}
 For any $x\in M$ and any neighborhood $\mathcal{U}_{x}$, there exists a uniformly normal neighborhood $\mathcal{V}_{x}$ such that $\mathcal{V}_{x}\subset \mathcal{U}_{x}$.
 \end{lemma}

 \section{Averaging on Riemannian Manifolds}
In this section we present the analysis of the averaging methods for nonlinear dynamical systems on Riemannian manifolds.
We derive the propagation equations for a single point under two different vector fields in order to bound the variation of the distance function between different state trajectories. 
\subsection{Closeness of Solutions}
Consider the following time-varying dynamical systems on $M$:
\EQ \label{1}&&\dot{x}(t)=f_{1}(x(t),t),\dot{y}(t)=f_{2}(y(t),t),\nnum\\&& x(t_{0})=y(t_{0})=x_{0}\in M, f_{1},f_{2}\in \mathfrak{X}(M\times \mathds{R}),\EN
where $\mathfrak{X}(M\times \mathds{R})$ is the space of smooth time-varying vector fields on $M$. 

\begin{theorem}[Closeness of Solutions]
\label{tp}
Consider the system of dynamical equations given by (\ref{1}) on the time interval $[t_{0},t_{1}]$. Then, 
\EQ &&d(\Phi_{f_{1}}(t,t_{0},x_{0}),\Phi_{f_{2}}(t,t_{0},x_{0}))\nnum\\&&\leq K(t_{1}-t_{0})\exp [C(t-t_{0})],t\in[t_{0},t_{1}],\EN
for some $K,C\geq 0$.
\end{theorem}
\begin{pf}
Consider a piecewise smooth path $\gamma(\tau)\in M,\tau\in[0,1]$  as follows (Theorem \ref{t1} guarantees the existence of $\gamma$):
\EQ \gamma:[0,1]\rightarrow M, \gamma(0)=\gamma(1)=x_{0}.\EN
Define a time and parameter varying vector field $X\in\mathfrak{X}(\mathds{R}\times\mathds{R}\times M)$ as
\EQ\label{moj} &&X(\tau,t,x)=f_{2}(x,t)+\tau(f_{1}(x,t)-f_{2}(x,t)),\nnum\\&& \tau\in[0,1],t\in[0,\infty)\subset\mathds{R},x\in M.\EN
It is clear that $X(0,t,x)=f_{2}(x,t)$, $X(1,t,x)=f_{1}(x,t)$, while  $X$ is smooth with respect to $\tau,t$ and $x$.
Hence, $\Phi_{X(0,t,x)}(t,t_{0},x_{0})=\Phi_{f_{2}(x,t)}(t,t_{0},x_{0})$ and  $\Phi_{X(1,t,x)}(t,t_{0},x_{0})=\Phi_{f_{1}(x,t)}(t,t_{0},x_{0})$.

An admissible family of curves, $\Gamma$, corresponding to $\Phi_{X(\tau,t,x)}(t,t_{0},\gamma(\tau))$ is given by (see Figure \ref{fff2}) 
\EQ\label{af} \Gamma:[0,1]\times \mathds{R}\rightarrow M,\hspace{.2cm} \Gamma(\tau,t)\doteq\Phi_{X(\tau,t,x)}(t,t_{0},\gamma(\tau)).\nnum\\\EN
\begin{figure}
\begin{center}
\hspace*{-.5cm}\includegraphics[scale=.25]{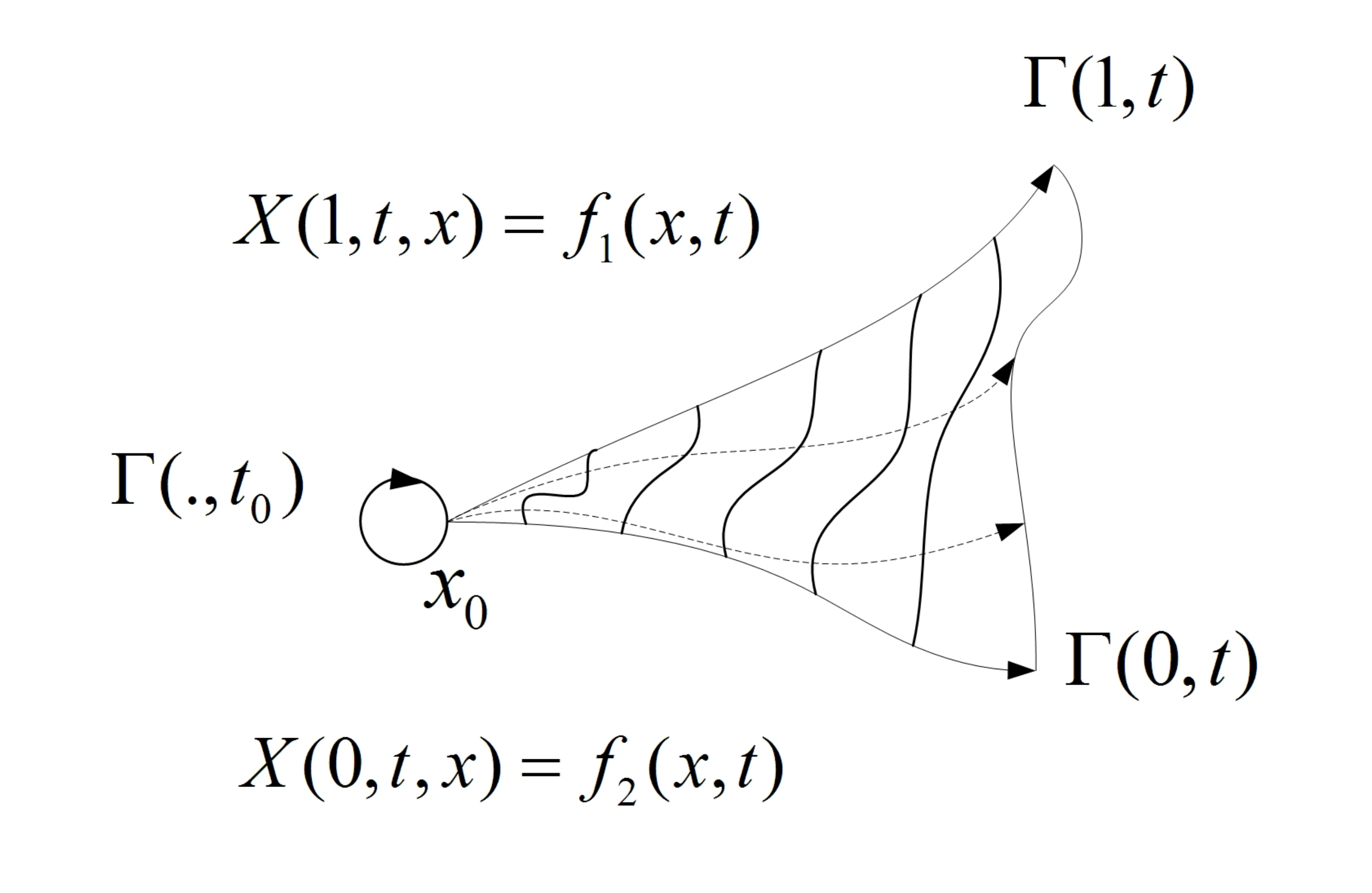}
 \caption{Admissible family of curves}
     \label{fff2}
      \end{center}
   \end{figure}  
Here we analyze the variation of $\ell(\Gamma(\cdot,t))$  with respect to $t$ where $\ell$ is the length function on $M$ where $\ell(\Gamma(t))\doteq \ell(\Gamma(\cdot,t))=\int^{1}_{0}||\partial_{\tau}\Gamma(\tau,t)||_{g}d\tau$  (for the sake of simplicity in our notation for this proof we drop the subscript $g$ from $||\cdot||_{g}$ in the following equations). In particular, note that
\EQ
\label{miri} \frac{d}{dt}\ell(\Gamma(t))&=&\frac{d}{dt}\int^{1}_{0}||\partial_{\tau}\Gamma(\tau,t)||d\tau\nnum\\&=&\frac{d}{dt}\int^{1}_{0}g(\nabla_{\frac{\partial}{\partial \tau}}\Gamma(\tau,t), \nabla_{\frac{\partial}{\partial \tau}}\Gamma(\tau,t))^{\frac{1}{2}}d\tau,\EN
where the second equality is implied by (\ref{levi}) (ii). Also note that
\EQ \label{kir2}\frac{\partial}{\partial t} g(\nabla_{\frac{\partial}{\partial \tau}}\Gamma(\tau,t), \nabla_{\frac{\partial}{\partial \tau}}\Gamma(\tau,t))^{\frac{1}{2}}&=&\frac{\partial}{\partial t}||\partial_{\tau}\Gamma(\tau,t)||\nnum\\&=&\frac{\frac{\partial}{\partial t}||\partial_{\tau}\Gamma(\tau,t)||2||\partial_{\tau}\Gamma(\tau,t)||}{2||\partial_{\tau}\Gamma(\tau,t)||}\nnum\\&=&\frac{\frac{\partial}{\partial t}||\partial_{\tau}\Gamma(\tau,t)||^{2}}{2||\partial_{\tau}\Gamma(\tau,t)||}.\EN
Hence, interchanging the order  of differentiation and integration in (\ref{miri}) and applying (\ref{kir2}), 
\EQ \frac{d}{dt}\ell(\Gamma(t)) &=&\int^{1}_{0}\frac{\frac{\partial}{\partial t}g(\nabla_{\frac{\partial}{\partial \tau}}\Gamma(\tau,t), \nabla_{\frac{\partial}{\partial \tau}}\Gamma(\tau,t))}{2||\nabla_{\frac{\partial}{\partial \tau}}\Gamma(\tau,t)||}d\tau\nnum\\&=&\int^{1}_{0}\frac{g(\nabla_{\frac{\partial}{\partial t}}\nabla_{\frac{\partial}{\partial \tau}}\Gamma(\tau,t), \nabla_{\frac{\partial}{\partial \tau}}\Gamma(\tau,t))}{||\nabla_{\frac{\partial}{\partial \tau}}\Gamma(\tau,t)||}d\tau\nnum\\&=&\int^{1}_{0}\frac{g(\nabla_{\frac{\partial}{\partial \tau}}\partial_{t}\Gamma(\tau,t), \nabla_{\frac{\partial}{\partial \tau}}\Gamma(\tau,t))}{||\nabla_{\frac{\partial}{\partial \tau}}\Gamma(\tau,t)||}d\tau\nnum\\&\leq &\int^{1}_{0}||\nabla_{\frac{\partial}{\partial \tau}}\partial_{t}\Gamma(\tau,t)||d\tau\nnum\\&=&\int^{1}_{0}||\nabla_{\frac{\partial}{\partial \tau}}X(\tau,t,\Gamma(\tau,t))||d\tau,\EN
where the second equality is by applying (\ref{kir4}) and (\ref{kir3}) together  and the inequality above is obtained by a direct application of Cauchy-Schwarz inequality. Hence, applying (\ref{moj}),
\EQ \frac{d}{dt}\ell(\Gamma(t))&\leq& \int^{1}_{0}||\nabla_{\frac{\partial}{\partial \tau}}\big[f_{2}(\Gamma(\tau,t),t)+\nnum\\&&\tau(f_{1}(\Gamma(\tau,t),t)-f_{2}(\Gamma(\tau,t),t))\big]||d\tau\nnum\\&=&\int^{1}_{0}||\nabla_{\frac{\partial}{\partial \tau}}f_{2}(\Gamma(\tau,t),t))+\nnum\\&&\nabla_{\frac{\partial}{\partial \tau}}\tau(f_{1}(\Gamma(\tau,t),t)-f_{2}(\Gamma(\tau,t),t))||d\tau\nnum\\&\leq&\int^{1}_{0}\big(||\nabla_{\frac{\partial}{\partial \tau}}f_{2}(\Gamma(\tau,t),t)||+\nnum\\&&||\nabla_{\frac{\partial}{\partial \tau}}\tau(f_{1}(\Gamma(\tau,t),t)-f_{2}(\Gamma(\tau,t),t))||\big)d\tau,\nnum\EN 
where the last inequality follows by an application of the triangle inequality. 
Employing (\ref{lev}),  
\EQ &&\nabla_{\frac{\partial}{\partial \tau}}\tau(f_{1}(\Gamma(\tau,t),t)-f_{2}(\Gamma(\tau,t),t))=\nnum\\&&+\tau\nabla_{\frac{\partial}{\partial \tau}}(f_{1}(\Gamma(\tau,t),t)-f_{2}(\Gamma(\tau,t),t))\nnum\\&&(f_{1}(\Gamma(\tau,t),t)-f_{2}(\Gamma(\tau,t),t)),\EN
so that another application of the triangle inequality yields
\EQ \frac{d}{dt}\ell(\Gamma(t))&\leq&\int^{1}_{0}||\nabla_{\frac{\partial}{\partial \tau}}f_{2}(\Gamma(\tau,t),t)||d\tau+\nnum\\&&\int^{1}_{0}||(f_{1}(\Gamma(\tau,t),t)-f_{2}(\Gamma(\tau,t),t))||d\tau+\nnum\\&&\int^{1}_{0}\tau||\nabla_{\frac{\partial}{\partial \tau}}(f_{1}(\Gamma(\tau,t),t)-f_{2}(\Gamma(\tau,t),t))||d\tau.\nnum\EN
Hence,
\EQ &&\ell(\Gamma(t))\leq \ell(\Gamma(t_{0}))+\int^{t}_{t_{0}}\int^{1}_{0}||\nabla_{\frac{\partial}{\partial \tau}}f_{2}(\Gamma(\tau,s),s)||d\tau ds+\nnum\\&&\int^{t}_{t_{0}}\int^{1}_{0}||(f_{1}(\Gamma(\tau,s),s)-f_{2}(\Gamma(\tau,s),s))||d\tau ds+\nnum\\&& \int^{t}_{t_{0}}\int^{1}_{0}\tau||\nabla_{\frac{\partial}{\partial \tau}}(f_{1}(\Gamma(\tau,s),s)-f_{2}(\Gamma(\tau,s),s))||d\tau ds.\nnum\EN
We note that since $\ell(\Gamma(t_{0}))=\ell(\gamma)$, we can choose the trivial path defined by $\gamma(\tau)=x_{0}, \tau\in[0,1]$. Hence, $\ell(\Gamma(t_{0}))=0$. Hence, without loss of generality, we have
\EQ\label{d} &&d(\Phi_{f_{1}}(t,t_{0},x_{0}),\Phi_{f_{2}}(t,t_{0},x_{0}))\leq \ell(\Gamma(t))\leq\nnum\\& &\int^{t}_{t_{0}}\int^{1}_{0}||\nabla_{\frac{\partial}{\partial \tau}}f_{2}(\Gamma(\tau,s),s)||d\tau ds\nnum\\&&+\int^{t}_{t_{0}}\int^{1}_{0}||(f_{1}(\Gamma(\tau,s),s)-f_{2}(\Gamma(\tau,s),s))||d\tau ds\nnum\\&&+ \int^{t}_{t_{0}}\int^{1}_{0}\tau||\nabla_{\frac{\partial}{\partial \tau}}(f_{1}(\Gamma(\tau,s),s)-f_{2}(\Gamma(\tau,s),s))||d\tau ds.\nnum\\\EN
Define
\EQ D_{\Gamma}\doteq\bigcup_{\tau\in[0,1], s\in[t_{0},t_{1}]}\Gamma(\tau,s)\subset M. \EN
Since $\Gamma$  is continuous on $[0,1],\times[t_{0},t_{1}]$, $D_{\Gamma}$ is  compact in the topology of $M$.  By our hypotheses $f_{1},f_{2}$ are smooth mappings and $\Gamma$ is continuous by construction. Therefore, $||f_{1}(x,t)-f_{2}(x,t)||$ attains its maximum on $D_{\Gamma}\times[t_{0},t_{1}],$ which is denoted by 
\EQ \label{tan}K_{\Gamma}\doteq\max_{(x,t)\in D_{\Gamma}\times[t_{0},t_{1}]}||f_{1}(x,t)-f_{2}(x,t)||.\EN
As is shown by (\ref{kos1}), the covariant differential of a vector field $X\in \mathfrak{X}(M)$, i.e. $\nabla X(x),\hspace{.2cm} x\in M$, is a linear operator as (see \cite{Lee3}, Chapter 4)
\EQ \nabla X(x):(T_{x}M,||.||)\rightarrow (T_{x}M,||.||).\EN
We denote the norm of this bounded linear operator by $||\nabla X(x)||$, so that 
\EQ \label{kos4}&&\exists C_{i}\in (0,\infty),\hspace{.2cm}s.t.\hspace{.2cm}\nnum\\&& C_{i}\doteq\sup_{x\in D_{\Gamma},t\in[t_{0},t_{1}]}||\nabla f_{i}(x,t)||,\hspace{.2cm}i=1,2.\EN 
It is shown in \cite{Lee3}, Lemma 4.2, that $\nabla_{X}Y(x),\hspace{.2cm} X,Y\in \mathfrak{X}(M),x\in M$, only depends on $X(x)\in T_{x}M$. Therefore,
\EQ \nabla_{\frac{\partial}{\partial \tau}}f_{i}(\Gamma(\tau,t),t)=\nabla_{\partial_{\tau}\Gamma(\tau,t)}f_{i}(\Gamma(\tau,t),t),\nnum\\ \tau\in[0,1],t\in[t_{0},t_{1}],i=1,2,\EN
since $\frac{\partial}{\partial \tau}|_{(\tau,t)}=\partial_{\tau}\Gamma(\tau,t)$.
Hence, 
\EQ \label{ying}||\nabla_{\frac{\partial}{\partial \tau}}f_{i}(\Gamma(\tau,t),t)||\leq C_{i}||\partial_{\tau}\Gamma(\tau,t)||.\EN
Applying (\ref{ying}) to (\ref{d}) yields
\EQ \label{jq}\ell(\Gamma(t))&\leq &\int^{t}_{t_{0}}\int^{1}_{0}C_{2}||\partial_{\tau}\Gamma(\tau,t)||d\tau ds\nnum\\&+&\int^{t}_{t_{0}}\int^{1}_{0}||(f_{1}(\Gamma(\tau,s),s)-f_{2}(\Gamma(\tau,s),s))||d\tau ds\nnum\\&+& \int^{t}_{t_{0}}\int^{1}_{0}\tau \big(C_{1}||\partial_{\tau}\Gamma(\tau,s)||+C_{2}||\partial_{\tau}\Gamma(\tau,s)||\big)d\tau ds\nnum\\&\leq &K_{\Gamma}(t-t_{0})+ \int^{t}_{t_{0}} \big(C_{1}+2C_{2}\big)\ell(\Gamma(s))ds\hspace{.2cm}\nnum\\&\leq&K_{\Gamma}(t_{1}-t_{0})+ \int^{t}_{t_{0}} \big(C_{1}+2C_{2}\big)\ell(\Gamma(s))ds,\EN
where to obtain the second inequality we employed  (\ref{tan}),(\ref{ying}) and $\int^{1}_{0}||\partial_{\tau}\Gamma(\tau,s)||d\tau=\ell(\Gamma(s))$.
The inequality (\ref{jq}) is in an appropriate form for an application of the Gronwall inequality, which yields
\EQ\label{amy} &&d(\Phi_{f_{1}}(t,t_{0},x_{0}),\Phi_{f_{2}}(t,t_{0},x_{0}))\leq \ell(\Gamma(t))\leq\nnum\\&& K_{\Gamma}(t_{1}-t_{0})\exp[(C_{1}+2C_{2})(t-t_{0})],\EN
with $K=K_{\Gamma}$ and $C=(C_{1}+2C_{2})$.\qed
\end{pf}
\subsection{Averaging on Perturbed Dynamical Systems}
Using the closeness of solutions Theorem \ref{tp}, averaging can be introduced for systems evolving on a
manifold $M$ (see Chapters 9 and 10 of \cite{Lewis} and \cite{Kha} respectively). In particular, we consider closeness
of solutions of two perturbed periodic systems of the form (\ref{1}), leading to the study of closeness of
solutions with respect to an averaged system. The resulting averaging Theorem is illustrated via a
subsequent application to a simple example. To this end, consider the following dynamical equations on a Riemannian manifold $M$:
\EQ \label{kos2}&&\dot{x}(t)=f^{\epsilon}_{1}(x,t)=\epsilon f_{1}(x,t),\hspace{.2cm} \nnum\\&&\dot{y}(t)=f^{\epsilon}_{2}(y,t)=\epsilon f_{2}(y,t),\nnum\\&&f_{1},f_{2}\in \mathfrak{X}(M\times \mathds{R}),\nnum\\&& x(t_{0})=y(t_{0})=x_{0}\in M \hspace{.2cm}x(t),y(t)\in M,0\leq\epsilon. \EN
The following lemma extends the closeness of solutions Theorem \ref{tp}  to perturbed dynamical systems on $M$. We note that the analyses presented in this paper can be extended to general non-periodic vector fields on Riemannian manifolds. In this case, the averaged vector fields are defined by averaging the nominal vector fields over an infinite time horizon, see \cite{nes4}.
\begin{lemma}
\label{lp}
Consider the dynamical systems of the form (\ref{kos2}) on $M$. Suppose there exists $\epsilon_{1}> 0$ such that the flows $\Phi_{\epsilon f_{i}}(\cdot,t_{0},x_{0}),\hspace{.2cm}i=1,2,$ exist on $[t_{0},t_{1}]$ for $\epsilon\in[0,\epsilon_{1}]$. Then for a time interval of order $O(1)$ and $\epsilon\in(0,\epsilon_{1}]$, we have
\EQ d(\Phi_{\epsilon f_{1}}(t,t_{0},x_{0}),\Phi_{\epsilon f_{2}}(t,t_{0}, x_{0}))=O(\epsilon),\hspace{.2cm}t\in[t_{0},t_{1}].\nnum\EN 
\end{lemma}
\begin{pf}
We define $\Gamma(\tau,t,\epsilon)$ as an admissible family of curves given by the flow of the vector field  $X(\tau,t,x,\epsilon)=\epsilon f_{2}(x,t)+\epsilon\tau(f_{1}(x,t)-f_{2}(x,t))\in T_{x}M,\hspace{.2cm} \tau\in[0,1],t\in[t_{0},t_{1}],x\in M$, such that  

\EQ\label{gg} \Gamma(\tau,t,\epsilon)\doteq\Gamma_{\epsilon}(\tau,t)=\Phi_{X(\tau,t,x,\epsilon)}(t,t_{0},\gamma(\tau))\in M, \hspace{.2cm}\nnum\\\gamma(\tau)=x_{0}, \tau\in[0,1],\EN
where $\Gamma_{\epsilon}(\tau,t)$ is of the same form as (\ref{af}).
By construction,  $\Gamma$ is continuous with respect to $(\tau,t)$. Employing the results of \cite{mar1}, it can be shown that $\Gamma$ is continuous with respect to $\epsilon$ as well. This yields compactness of $\hat{D}_{\Gamma}$, where   
\EQ \hat{D}_{\Gamma}\doteq\bigcup_{\tau\in[0,1], t\in[t_{0},t_{1}],\epsilon\in[0,\epsilon_{1}]}\Gamma(\tau,t,\epsilon).\EN
We then modify $K_{\Gamma}$ and $C_{i},\hspace{.2cm}i=1,2$, as per (\ref{tan}) and (\ref{kos4}), to define
\EQ  &&\hat{K}_{\Gamma}\doteq\sup_{\hat{D}_{\Gamma}\times[t_{0},t_{1}]}||f_{1}(x,t)-f_{2}(x,t)||,\nnum\\&& \hat{C}_{i}\doteq\sup_{\hat{D}_{\Gamma}\times[t_{0},t_{1}]}||\nabla f_{i}(x,t)||,\hspace{.2cm}i=1,2.\EN
Applying Theorem \ref{tp} then yields 
\EQ &&d(\Phi_{f_{1}}(t,t_{0},x_{0}),\Phi_{f_{2}}(t,t_{0},x_{0}))\leq\nnum\\&&  \epsilon\hat{K}_{\Gamma}(t_{1}-t_{0})\exp[\epsilon_{1}(\hat{C}_{1}+2\hat{C}_{2})(t-t_{0})]=O(\epsilon),\nnum\EN
which completes the proof.\qed
 \end{pf}

 Let us consider a perturbed system as
\EQ \label{pp}\dot{x}(t)=\epsilon f(x(t),t), \hspace{.2cm}f\in \mathfrak{X}(M\times \mathds{R}),x_{0}\in M,\hspace{.2cm}\epsilon\geq 0,\nnum\\\EN
where $f$ is periodic in $t$ with the period $T$, i.e. $f(x,t)=f(x,t+T)$. Such a system is referred to as \textit{$T$-periodic}. The averaged vector field $\hat{f}$ is given by 
\EQ \label{ppp}\hat{f}(x)\doteq\frac{1}{T}\int^{T}_{0}f(x,s)ds,\EN
where the average dynamical system is locally given by $\dot{x}(t)=\epsilon \hat{f}(x(t))$.
The following theorem is the first order averaging theorem for periodic dynamical systems on compact Riemannian manifolds.
\begin{theorem}[Averaging Theorem]
\label{ta}
For a smooth $n$ dimensional compact Riemannian manifold $M$, let $f\in \mathfrak{X}(M\times \mathds{R})$ be a $T$-periodic smooth vector field. Then, for any given $t_{1}\in [t_{0},\infty)$, such that $t_{1}-t_{0}=O(\frac{1}{\epsilon}),\hspace{.2cm}\epsilon\in(0,\epsilon_{1}]$ for some  $0<\epsilon_{1}$,
\EQ d(\Phi_{\epsilon f}(t,t_{0},x_{0}),\Phi_{\epsilon \hat{f}}(t,t_{0},x_{0}))=O(\epsilon).\EN
\end{theorem}
In order to prove Theorem \ref{ta}, we employ the notion of pullbacks of vector fields along diffeomorphisms on 
$M$. Let $X,Y\in \mathfrak{X}(M\times \mathds{R})$ be smooth time-varying vector fields on $M$, where it may be shown that $\Phi_{Y}(t,t_{0},.):M\rightarrow M$ is a local diffeomorphism (see \cite{mar1}). Define
\EQ\label{pu} &&\Phi_{Y}^{(t,t_{0})^{*}}:T_{\Phi_{Y}(t,t_{0},x_{0})}M\rightarrow T_{x_{0}}M,\nnum\\&& \Phi_{Y}^{(t,t_{0})^{*}}X(x_{0},s)\doteq T\Phi_{Y}^{(t,t_{0})^{-1}}X(\Phi_{Y}(t,t_{0}, x_{0}),s),\nnum\\&& t,s\in\mathds{R}, X\in \mathfrak{X}(M\times \mathds{R}),\EN
where $T\Phi_{Y}^{(t,t_{0})^{-1}}$ is the pushforward of $\Phi_{Y}^{-1}(t,t_{0},.):M\rightarrow M$ defined in the standard framework of differential geometry (see \cite{Lee2}, Chapter 3).
We have the following lemma for the variation of smoothly varying vector fields with respect to a parameter variable.
\begin{lemma}[\hspace{-.02cm}\cite{Agra}, Page 40, \cite{Lewis}, Page 451]
\label{ll1}
Consider a smooth vector field $Y\in \mathfrak{X}(\mathds{R}\times M)$ with the associated flow $\Phi_{Y}(t,t_{0},\cdot):M\rightarrow M$. Then,
\EQ&& \frac{\partial}{\partial\lambda}\Phi_{Y(\lambda,x(t))}(t,t_{0},x_{0})=T_{x_{0}}\Phi^{(t,t_{0})}_{Y(\lambda,x(t))}\times\nnum\\&&\int^{t}_{t_{0}}\Phi^{(s,t_{0})^{*}}_{Y(\lambda,x(s))}\frac{\partial }{\partial \lambda}Y\left(\lambda,\Phi_{Y(\lambda,x(s))}(s,t_{0},x_{0})\right)ds=\nnum\\&&\int^{t}_{t_{0}}(\Phi^{-1})^{(t,s)^{*}}_{Y(\lambda,x(s))}\frac{\partial }{\partial \lambda}Y(\lambda,\Phi_{Y(\lambda,x(s))}(s,t_{0},x_{0}))ds\nnum\\&& \in T_{\Phi_{Y(\lambda,x(t))}^{(t,t_{0})}(x_{0})}M.\EN
\hspace*{8cm}\qed\end{lemma}
The proof of Theorem 4 follows via the methodology of \cite{Lewis} for dynamical systems evolving
on $\mathds{R}^n$, and an extension of the results of \cite{Kha}, Theorem 10.4.

\begin{pf}(Theorem \ref{ta}) 
Define the smooth parameter varying vector field $Y(\lambda,x)\doteq\int^{\lambda}_{0}\big(f(x,s)-\hat{f}(x)\big)ds,\hspace{.2cm}0\leq \lambda,$ then by Lemma \ref{ll1} we have
\EQ &&\frac{\partial}{\partial\lambda}\Phi_{\epsilon Y(\lambda,x(t))}(t,t_{0},x_{0})=\nnum\\&&\epsilon\int^{t}_{t_{0}}(\Phi^{-1})^{(t,s)^{*}}_{\epsilon Y(\lambda,x(s))}\frac{\partial }{\partial \lambda}Y(\lambda,\Phi_{Y(\lambda,x(s))}(s,t_{0},x_{0}))ds=\nnum\\&&\epsilon\int^{t}_{t_{0}}(\Phi^{-1})^{(t,s)^{*}}_{\epsilon Y(\lambda,x(s))}\big(f(x(s),\lambda)-\hat{f}(x(s))\big)ds.\EN

For a given initial condition $y_{0}\in M$, define a perturbed curve $y:\mathds{R}\times \mathds{R}\rightarrow M $ by
\EQ y(\lambda,\tau)\doteq\Phi_{\epsilon Y(\lambda,y)}(\tau,0,y_{0}),\hspace{.2cm}\nnum\\ y_{0}\in M,\tau\in[0,1],\lambda\in\mathds{R}_{\geq 0}.\EN

Since $Y(\lambda,x)$ is smooth with respect to both $x$ and $\lambda$, then $\Phi_{\epsilon Y(\lambda,y(\lambda,\tau))}(1,0,y_{0})$ has the same degree of regularity with respect to $\lambda$ (see \cite{mar1, Lewis}, Page 450). Note that the existence of $\Phi_{\epsilon Y(\lambda,y(\lambda,\tau))}(1,0,y_{0})$ is guaranteed by the compactness of 
$M$. Define
\EQ \label{b}D_{\Phi,\epsilon}\doteq\bigcup_{\tau\in[0,1]}\Phi_{\epsilon Y(\lambda,y)}(\tau,0,y_{0})\subset M, \hspace{.2cm}\lambda\in\mathds{R}_{\geq 0}.\EN
Now we show that
\EQ d(\Phi_{\epsilon Y(t,x(t))}(1,0,y_{0}),y_{0})=O(\epsilon),\hspace{.2cm}t\in[t_{0},\infty).\EN
By the definition of the length function in (\ref{l}),\\  $d(\Phi_{\epsilon Y(t,x(t))}(t,t_{0},y_{0}),y_{0})\leq \ell(\Phi_{\epsilon Y(t,x(t))}(t,t_{0},y_{0})),$ therefore
\EQ &&d(\Phi_{\epsilon Y(\lambda,y(\lambda,\tau))}(1,0,y_{0}),y_{0})\leq \nnum\\&&\ell(\Phi_{\epsilon Y(\lambda,y(\lambda,\tau))}(1,0,y_{0}))\leq \epsilon\int^{1}_{0}||Y(\lambda,y(\lambda,\tau))||d\tau.\nnum\EN
Periodicity of $Y$ with respect to $\lambda$, boundedness of $y(\lambda,\tau),\hspace{.2cm}\lambda\in[0,T],$ in the sense of prempactness of $D_{\Phi,\epsilon}$ (i.e. $D_{\Phi,\epsilon}$ is contained in a compact set $M$) in (\ref{b}) and smoothness of $Y$ with respect to $y$ together yield $d(\Phi_{\epsilon Y(\lambda,y(\lambda,\tau))}(1,0,y_{0}),y_{0})=O(\epsilon)$.

In order to obtain the statement of the theorem it is sufficient to prove that
\EQ d(\Phi_{\epsilon f(x,t)}(t,t_{0},x_{0}),\Phi_{\epsilon Y(t,x)}^{(1,0)}\circ \Phi_{\epsilon \hat{f}(x)}(t,t_{0},x_{0}))=O(\epsilon),\nnum\EN
since by the triangle inequality,
\EQ &&d(\Phi_{\epsilon f(x,t)}(t,t_{0},x_{0}), \Phi_{\epsilon \hat{f}(x)}(t,t_{0},x_{0}))\leq\nnum\\&& d(\Phi_{\epsilon f(x,t)}(t,t_{0},x_{0}),\Phi_{\epsilon Y(t,x)}^{(1,0)}\circ \Phi_{\epsilon \hat{f}(x)}(t,t_{0},x_{0}))+\nnum\\&& d(\Phi_{\epsilon Y(t,x)}^{(1,0)}\circ \Phi_{\epsilon \hat{f}(x)}(t,t_{0},x_{0}),\Phi_{\epsilon \hat{f}(x)}(t,t_{0},x_{0})),\EN
where \\$\Phi_{\epsilon Y(t,x)}^{(1,0)}\circ \Phi_{\epsilon \hat{f}(x)}(t,t_{0},x_{0})\doteq\Phi_{\epsilon Y(t,x)}(1,0, \Phi_{\epsilon \hat{f}(x)}(t,t_{0},x_{0}))$.
Here we compute the tangent vector field of $y(t)=\Phi_{\epsilon Y(t,x)}^{(1,0)}\circ \Phi_{\epsilon \hat{f}(x)}(t,t_{0},x_{0})\in M$. The derivative of $\Phi_{\epsilon Y(t,x)}^{(1,0)}\circ \Phi_{\epsilon \hat{f}(x)}(t,t_{0},x_{0})$ with respect to time can be computed via the chain rule as follows:
\EQ \label{rr}&&\dot{y}(t)=T_{\Phi_{\epsilon \hat{f}}(t,t_{0},x_{0})}\Phi_{\epsilon Y(t,y)}^{(1,0)}\Big(\epsilon \hat{f}(\Phi_{\epsilon \hat{f}}(t,t_{0},x_{0}))\Big)+\nnum\\&&\frac{\partial}{\partial t}\Phi_{\epsilon Y(t,y)}^{(1,0)}\big(\Phi_{\epsilon Y(t,y)}^{(1,0)}\circ \Phi_{\epsilon \hat{f}(x)}(t,t_{0},x_{0})\big)=(\Phi^{-1})_{\epsilon Y(t,y)}^{(1,0)^{*}}\nnum\\&&\Big(\epsilon \hat{f}(\Phi_{\epsilon \hat{f}}(t,t_{0},x_{0}))\Big)+\epsilon\int^{1}_{0}(\Phi^{-1})^{(1,s)^{*}}_{\epsilon Y(t,y(s))}\nnum\\&&\big(f(y(s),t)-\hat{f}(y(s))\big)ds,\EN 
where the second equality is established by the definition of pullbacks in (\ref{pu}) and the equation of parameter variation of flows given by Lemma \ref{ll1}.
In a compact form, (\ref{rr}) is written as
\EQ \dot{y}(t)&=&\epsilon\Big[(\Phi^{-1})_{\epsilon Y(t,y)}^{(1,0)^{*}} \hat{f}+\int^{1}_{0}(\Phi^{-1})^{(1,s)^{*}}_{\epsilon Y(t,y(s))}\nnum\\&&\big(f-\hat{f}\big)ds\Big]\circ y(t)\doteq\epsilon G(\epsilon,t,y(t)).\EN
where $G(\epsilon,t,y)\in T_{y(t)}M$. Since the vector fields $f,\hat{f}$ are both smooth, the construction above implies that $G$ is  smooth with respect to $\epsilon$. One can see that by setting $\epsilon=0$, the nominal vector field $f$ is retrieved from $G$, i.e. $G(0,t,x)=f(x,t)$. This is due to the fact that at $\epsilon=0$, $(\Phi^{-1})_{\epsilon Y(t,x)}^{(1,0)^{*}}=I$ and the state trajectory $y(.)$ will be independent of $s$ in the integral term of $G$ (for an identically zero vector field, the state trajectory does not evolve away from its initial state). By applying the Taylor expansion with remainder we have
\EQ G(\epsilon,t,x)=f(x,t)+\epsilon h(x,\zeta,t),\EN
where $h(x,\zeta,t)=\frac{\partial }{\partial \epsilon}G(\epsilon,t,x)|_{\epsilon=\zeta}$ and $\zeta\in[0,\epsilon]$ and both $G$ and $h$ are $T$ periodic.
Now let us explore the state variation along the following dynamical equations:
\EQ &&\dot{x}(t)=\epsilon f(x(t),t),\hspace{2.2cm}x(t_{0})=x_{0},\nnum\\&&\dot{y}(t)=\epsilon f(y(t),t)+\epsilon^{2} h(y,\zeta,t),\hspace{.2cm}y(t_{0})=x_{0}.\EN
We note that $G$ is a smooth vector field on a compact Riemannian manifold $M$. Therefore, employing the results of the Escape Lemma (see \cite{Lee2}, Lemma 17.10) yields completeness of the flow of $G$ on $M$.
Following the results  of Theorem \ref{tp} and (\ref{jq}),

\EQ\label{ggg} &&d(\Phi_{\epsilon f(x,t)}(t,t_{0},x_{0}),\Phi_{\epsilon Y(t,x)}^{(1,0)}\circ \Phi_{\epsilon \hat{f}(x)}(t,t_{0},x_{0}))\leq\nnum\\&& \epsilon^{2}K_{\Gamma,h}(t-t_{0})\exp[\epsilon(C+\epsilon\hat{C})(t-t_{0})],\nnum\\\EN
where there exist $0<K_{\Gamma,h},C,\hat{C}<\infty,$ such that
\EQ\label{para} && K_{\Gamma,h}\doteq\sup_{(x,t)\in M\times[t_{0},t_{0}+T]}||h(x,\zeta,t)||,\nnum\\&& C\doteq\sup_{ (x,t)\in M\times[t_{0},t_{0}+T]}||\nabla f(x,t)||,\nnum\\&&\hspace{.2cm}\hat{C}\doteq\sup_{(x,t)\in M\times[t_{0},t_{0}+T]}||\nabla h(x,\zeta,t)||.\EN
The parameters $K_{\Gamma,h}, C,$ and $\hat{C}$ are all invariant with respect to $x$, since  
\EQ \hat{D}^{\infty}_{\Gamma}\doteq\bigcup_{\tau\in[0,1], t\in[t_{0},\infty],\epsilon\in[0,\epsilon_{1}]}\Gamma(\tau,t,\epsilon)\subset M,\EN
where $\Gamma$ is defined by (\ref{gg}). Also $f$ and $h$ are both $T$ periodic therefore the maximization in (\ref{para}) is taken on $t\in[t_{0},t_{0}+T]$.
\\\\Obviously for $t-t_{0}=O(\frac{1}{\epsilon})$ we have 
\EQ d(\Phi_{\epsilon f(x,t)}(t,t_{0},x_{0}),\Phi_{\epsilon Y(t,x)}^{(1,0)}\circ \Phi_{\epsilon \hat{f}(x)}(t,t_{0},x_{0}))=O(\epsilon),\nnum\EN
which completes the proof.\qed
\end{pf}
   \subsection{An example of averaging on $SO(3)$}
   \label{ex1}
In this section we present an example on $SO(3)$ which is a compact Lie group, see \cite{Lee2}.
The Lie algebra $\mathcal{L}$ of a Lie group $G$ is the tangent space at the identity element $e$ with the associated Lie bracket defined on the tangent space of $G$, i.e.  $\mathcal{L}=T_{e}G$. 
A vector field $X$ on $G$ is called \textit{left invariant } if
\EQ \forall g_{1},g_{2}\in G,\quad X(g_{1}\star g_{2})=TL_{g_{1}}X(g_{2}),\EN
where $L_{g_{1}}:G\rightarrow G,\hspace{.2cm} L_{g_{1}}(h)=g_{1}\star h,\hspace{.2cm} TL_{g_{1}}:T_{g_{2}}G\rightarrow T_{g_{1}\star g_{2}}G$ which immediately imply $X(g_{1}\star e)=X(g_{1})=TL_{g_{1}}X(e)$.
\begin{figure}
        \begin{center}
\hspace*{-.25cm}\includegraphics[scale=.25]{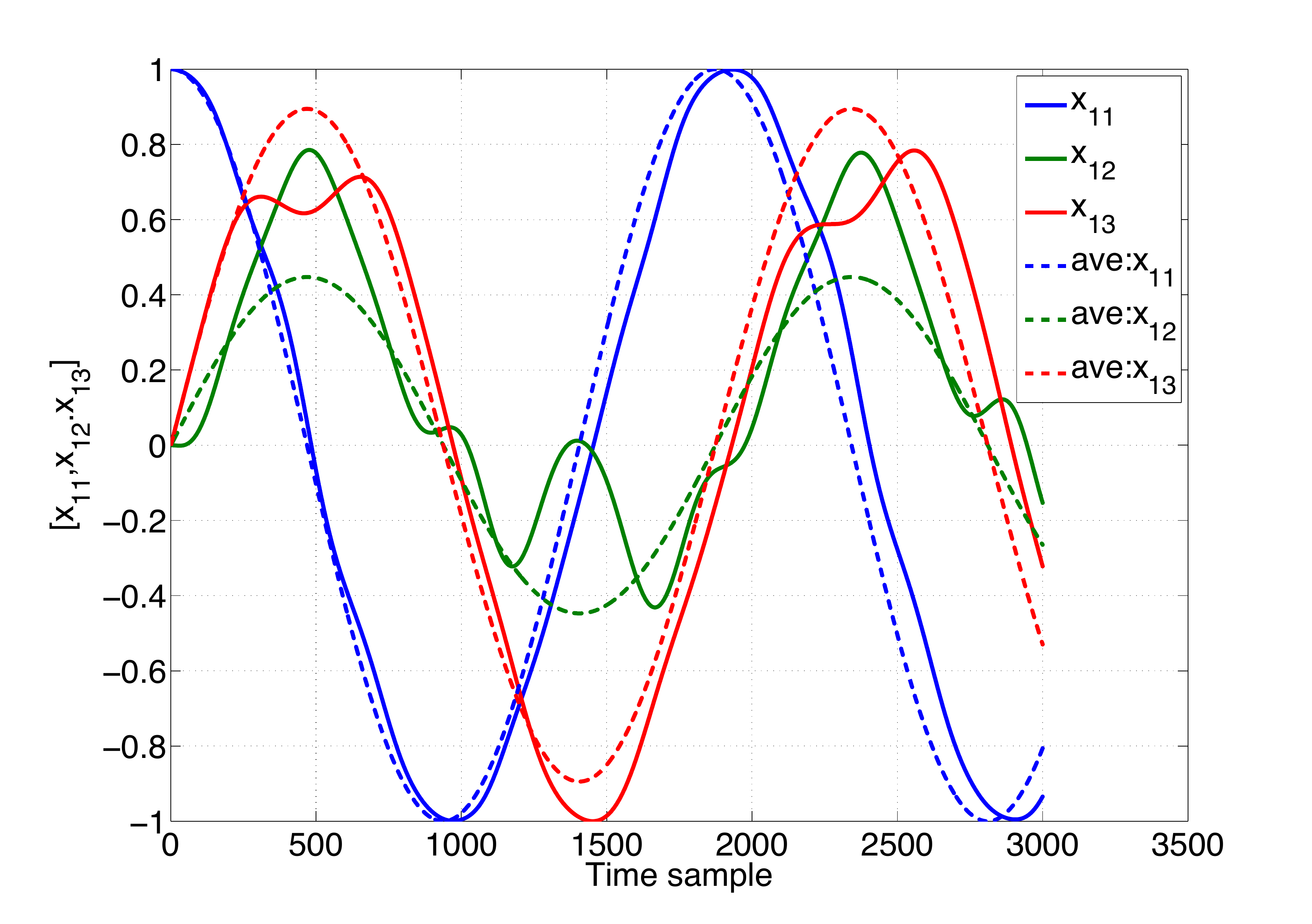}
 \caption{State trajectories on $SO(3)$(nominal system: solid line, average system: dashed line), $\epsilon=0.5$ (see Section \ref{ex1})}
     \label{ff11}
      \end{center}
   \end{figure}  
   \begin{figure}
\begin{center}
\hspace*{-.25cm}\includegraphics[scale=.25]{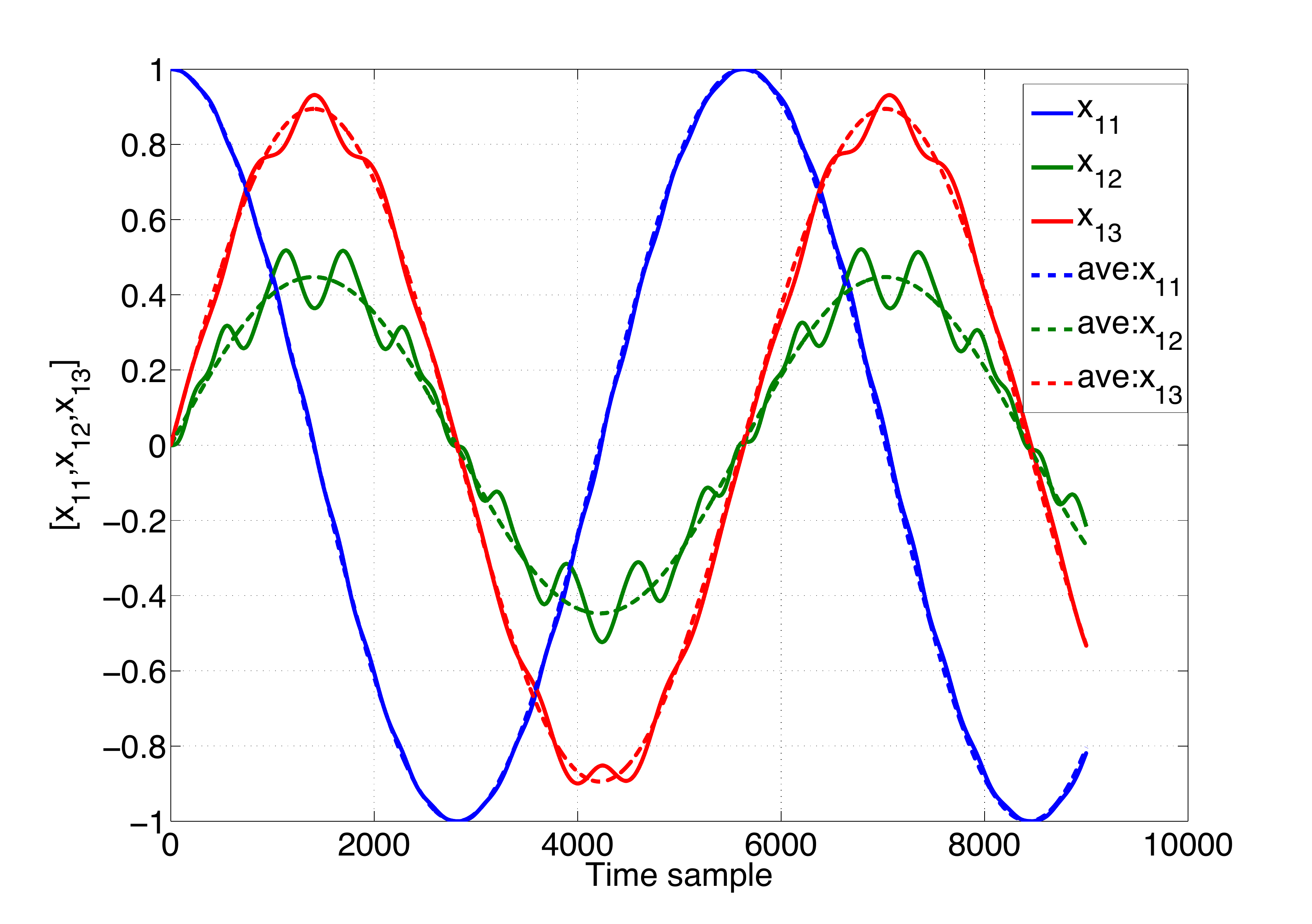}
 \caption{State trajectories on $SO(3)$ (nominal system: solid line, average system: dashed line), $\epsilon=0.1$(see Section \ref{ex1})}
     \label{ff22}
      \end{center}
   \end{figure}  
        We recall that $SO(3)$ is the rotation group in $\mathds{R}^{3}$ given by 
\EQ SO(3)= \big \{x\in GL(3)\hspace{.1cm}| \quad x\,\cdot\, x^{T}=I,\hspace{.1cm} det(x)=1\big\},\EN 
where $GL(n)$ is the set of nonsingular $n\times n$ matrices.
The Lie algebra of $SO(3)$ which is denoted by $so(3)$ is given by (see \cite{Varad})
\EQ so(3)=\big\{     X\in M(3)\hspace{.1cm}|\quad X+X^{T}=0\big\},\EN
where $M(n)$ is the space of all $n\times n$ matrices. The Lie group operation $\star$ is given by the matrix multiplication and consequently $TL_{g_{2}}$ is also given by the matrix multiplication $g_{2}X,\hspace{.2cm} X\in T_{g_{1}}G$.

A left invariant dynamical system on $SO(3)$ is given by (for the definition of left invariant dynamical systems see \cite{Lewis})
\EQ \dot{x}(t)=x(t)X(t),\quad x(0)=x_{0},\hspace{.2cm} X(t)\in so(3).\EN
The Lie algebra bilinear operator is defined as the commuter of matrices, i.e. $[X,Y]=XY-YX,\quad X,Y\in so(3).$ 
A controlled left invariant system on $SO(3)$ is then defined by
\EQ &&\hspace{-.2cm}\begin{array}{ll}\dot{x}(t)=x(t)\left( \begin{array}{ll} 0 \quad \hspace{.5cm}u_{1}(t) \quad u_{3}(t)\\-u_{1}(t)\quad 0\quad\hspace{.3cm} u_{2}(t)\\-u_{3}(t)\quad \hspace{-.2cm}-u_{2}(t)\quad 0
        \end{array}\right), \hspace{.2cm}\nnum\\x(t)\in SO(3), (u_{1}(t),u_{2}(t),u_{3}(t))\in \mathds{R}^{3}.\end{array}\EN
        The Lie algebra $so(3)$ is spanned by $e_{1}=\left( \begin{array}{ll} 0 \quad\hspace{.3cm} 1 \quad 0\\-1\quad 0\quad 0\\0\quad\hspace{.3cm} 0\quad 0
        \end{array}\right),\nnum\\ e_{2}=\left( \begin{array}{ll} 0 \quad\hspace{.3cm} 0 \quad 0\\0\quad\hspace{.3cm} 0\quad 1\\0\quad -1\quad 0
        \end{array}\right) \mbox{and}\hspace{.2cm}e_{3}=\left( \begin{array}{ll} 0 \quad\hspace{.3cm} 0 \quad 1\\0\quad\hspace{.3cm} 0\quad 0\\-1\quad 0\quad 0
        \end{array}\right)$.
        Consider the following perturbed left invariant dynamical system on $SO(3)$:
        
         \EQ &&\hspace{-.2cm}\dot{x}(t)=\epsilon x(t)\left( \begin{array}{ll} 0 \quad \hspace{.7cm}\sin^{2}(t) \quad 1\\-\sin^{2}(t)\quad 0\quad\hspace{.3cm} \cos(t)\\-1 \quad \hspace{.2cm}-\cos(t)\quad 0\end{array}\right).\EN
         The average dynamical system is given by $\dot{x}(t)=\epsilon x(t)\left( \begin{array}{ll} 0 \quad \hspace{.3cm}\frac{1}{2} \quad 1\\-\frac{1}{2}\quad 0\quad\hspace{0cm} 0\\-1 \quad \hspace{0cm}0\quad 0\end{array}\right).$
         Figures \ref{ff11} and \ref{ff22} show the closeness of solutions for the nominal and averaged systems above for $\epsilon=.5,$ and $.1$ respectively for $t\in[0,20]$ and $t\in[0,100]$ as expected by the results of Theorem \ref{ta}. 

\section{Infinite horizon averaging on Riemannian manifolds}
\label{s4}
Closeness of solutions on a finite time horizon may be extended to the infinite horizon limit via the
incorporation of appropriate stability properties, yielding averaging results for systems evolving on (not
necessarily compact) Riemannan manifolds. To this end, it is useful to state a number of standard
stability properties defined with respect to such manifolds.
\begin{definition}
For the dynamical system $\dot{x}(t)=f(x), \hspace{.2cm}f:M\rightarrow TM$,  $\bar{x}\in M$ is an equilibrium if 
\EQ \Phi_{f}(t,t_{0},\bar{x})=\bar{x},\hspace{.2cm}t\in[t_{0},\infty),\EN
where $\Phi_{f}$ is the flow of $f$.
\qed\end{definition}
\begin{definition}[\hspace{-.02cm}\cite{Lewis, Kha,forni,Ang}]
For the dynamical system $\dot{x}(t)=f(x), \hspace{.2cm}f:M\rightarrow TM$, an equilibrium $\bar{x}\in M$ is\\ 
(i): \textit{Lyapunov stable} if for any $t_0\in\mathds{R}$ and any neighborhood $\mathcal{U}_{\bar{x}}$ of $\bar{x}$, there exits a neighborhood $\mathcal{W}_{\bar{x}}$ of $\bar{x}$, such that 
\EQ x(t_{0})\in \mathcal{W}_{\bar{x}}\Rightarrow \Phi_{f}(t,t_{0},x(t_{0}))\in \mathcal{U}_{\bar{x}},\hspace{.2cm}t\in[t_{0},\infty).\EN\\
(ii): locally asymptotically stable if it is Lyapunov stable and for any $t_0\in\mathds{R},$ there exits $\mathcal{U}_{\bar{x}},$ such that
\EQ \forall x(t_{0})\in \mathcal{U}_{\bar{x}},\hspace{.2cm}\lim_{t\rightarrow \infty}\Phi_{f}(t,t_{0},x(t_{0}))=\bar{x}. \EN\\
(iii) globally asymptotically stable if it is Lyapunov stable and for any $t_0\in\mathds{R}$,  
\EQ \forall x(t_{0})\in M,\hspace{.2cm}\lim_{t\rightarrow \infty}\Phi_{f}(t,t_{0},x(t_{0}))=\bar{x}. \EN
(iv): locally exponentially stable if it is locally asymptotically stable and for any $t_0\in\mathds{R}$, there exists $\mathcal{U}_{\bar{x}},$ such that 
\EQ d(\Phi_{f}(t,t_{0},x(t_{0})),\bar{x})\leq kd(x(t_{0}),\bar{x})\exp(-\lambda(t-t_{0})),\nnum\\ k,\lambda\in \mathds{R}_{>0},x(t_{0})\in \mathcal{U}_{\bar{x}}.\nnum\\\EN \hspace*{7.8cm}\qed
\end{definition}
We note that the convergence on $M$ is defined in the topology induced by the metric $d$ which is same as the original topology of $M$ by Theorem \ref{t2}.  

\begin{definition}[\hspace{-.02cm}\cite{Lewis, Kha}]
A function $v:M\rightarrow \mathds{R}$ is locally positive-definite (positive-semidefinite) around $\bar{x}$ if $v(\bar{x})=0$ and there exists a neighborhood $\mathcal{U}_{\bar{x}}\subset M$ such that for all
$ x\in \mathcal{U}_{\bar{x}}-\{\bar{x}\},\hspace{.2cm} 0<v(x) \hspace{.2cm}(\mbox{respectively}\hspace{.2cm} 0\leq v(x)).$\qed
\end{definition}
Given a smooth function $v:M\rightarrow \mathds{R}$, the Lie derivative of $v$ along a vector field $f$ is defined by
\EQ \mathfrak{L}_{f}v\doteq dv(f),\EN
where $dv:TM\rightarrow \mathds{R}$ is the differential form of $v$, locally given by (see \cite{Lee2})
\EQ dv=\sum^{n}_{i=1}\frac{\partial v}{\partial x_{i}}dx_{i},\EN
where $n\doteq dim(M)$.
   
\begin{definition}
A smooth function $v:M\rightarrow \mathds{R}$ is a Lyapunov function for the vector field $f$, if $v$ is locally positive definite around the equilibrium $\bar{x}$ and $\mathfrak{L}_{f}v$ is locally negative-definite.  \qed
\end{definition}
\begin{definition}
The sublevel set $\mathcal{N}_{b}$ of a positive semidefinite function $v:M\rightarrow\mathds{R}$ is defined as $\mathcal{N}_{b}\doteq\{x\in M, \hspace{.2cm} v(x)\leq b\}$. By $\mathcal{N}_{b}(\bar{x})$ we denote the connected sublevel set of $M$ containing $\bar{x}\in M$. \hspace*{7.8cm}\qed
\end{definition}

The following lemma shows that there exists a connected compact neighborhood of an equilibrium point of a dynamical system on a Riemannian manifold.
\begin{lemma}[\hspace{-.02cm}\cite{Lewis}]
\label{lc}
Let $\bar{x}\in M$ be an equilibrium of $\dot{x}=f(x(t)),\hspace{.2cm}x(t)\in M$ and 
$v$ be a Lyapunov function on a neighborhood of $\bar{x}$. Then, for any neighborhood $\mathcal{U}_{\bar{x}}$ of $\bar{x}$, there exists $b\in\mathds{R}_{> 0}$ such that $\mathcal{N}_{b}(\bar{x})$ is compact, $\bar{x}\in int(\mathcal{N}_{b}(\bar{x}))$ and $\mathcal{N}_{b}(\bar{x})\subset \mathcal{U}_{\bar{x}}$.\hspace*{7.3cm}\qed  
\end{lemma}

\begin{theorem}
\label{taa}
For a smooth $n$ dimensional Riemannian manifold $M$, let $f\in \mathfrak{X}(M\times \mathds{R})$ be a $T$-periodic smooth vector field and assume  the nominal and averaged vector fields are both complete for $\epsilon\in(0,\epsilon_{1}],\hspace{.2cm}0< \epsilon_{1}$. Suppose the averaged dynamical system has a locally exponentially stable equilibrium $\bar{x}\in M$ such that there exists a Lyapunov function $v:M\rightarrow \mathds{R}$  where $\mathfrak{L}_{\hat{f}}v$ is locally negative-definite around $\bar{x}$. Then, there exists a neighborhood $\mathcal{N}_{\bar{x}}$ and $\hat{\epsilon}\leq\epsilon_{1}$ such that   
\EQ d(\Phi_{\epsilon f}(t,t_{0},x_{0}),\Phi_{\epsilon \hat{f}}(t,t_{0},x_{0}))=O(\epsilon),\hspace{.2cm}\nnum\\\epsilon\in(0,\hat{\epsilon}],x_{0}\in \mathcal{N}_{\bar{x}},\hspace{.2cm} t\in[t_{0},\infty),\EN
where $\hat f$ is the averaged vector field (\ref{ppp}).
\end{theorem}
 \begin{pf}
 First we note that the existence of a Lyapunov function $v:M\rightarrow \mathds{R}$ around $\bar{x}$, where $\mathfrak{L}_{\hat{f}}v$ is locally negative-definite around $\bar{x}$, guarantees that $\bar{x}$ is locally asymptotically stable (see \cite{Lewis}, Theorem 6.14). 
In order to analyze the dynamical system (\ref{pp}) on $[t_{0},\infty)$, we subtract the nominal vector field  from the averaged vector field and integrate, yielding 
\EQ Z(\lambda,x)\doteq\int^{\lambda}_{0}(\hat{f}(x)-f(x,\tau))d\tau,\hspace{.2cm}x\in M,\lambda\in\mathds{R}_{\geq 0}.\nnum\\\EN
Now consider a composition of flows on $M$ given by:
\EQ z(t)=\Phi^{(1,0)}_{\epsilon Z(y,t)}\circ \Phi_{\epsilon f(x,t)}(t,t_{0},x_{0}).\EN
Similar to the proof of Theorem \ref{ta}, the tangent vector field of $z$ is computed by
\EQ \label{rrr}\dot{z}(t)&=&T_{\Phi_{\epsilon f}(t,t_{0},x_{0})}\Phi_{\epsilon Z(t,z)}^{(1,0)}\Big(\epsilon f(\Phi_{\epsilon f}(t,t_{0},x_{0}),t)\Big)\nnum\\&&+\frac{\partial}{\partial t}\Phi_{\epsilon Z(t,z)}^{(1,0)}\big(\Phi_{\epsilon Z(t,z)}^{(1,0)}\circ \Phi_{\epsilon f(x,t)}(t,t_{0},x_{0})\big)\nnum\\&=&(\Phi^{-1})_{\epsilon Z(t,z)}^{(1,0)^{*}}\Big(\epsilon f(\Phi_{\epsilon f}(t,t_{0},x_{0}),t)\Big)\nnum\\&&+\epsilon\int^{1}_{0}(\Phi^{-1})^{(1,s)^{*}}_{\epsilon Z(t,z(s))}\big(\hat{f}(z(s))-f(z(s),t)\big)ds,\nnum\\\EN 
or equivalently 
\EQ \label{kk}\dot{z}(t)&=&\epsilon\Big[(\Phi^{-1})_{\epsilon Z(t,z)}^{(1,0)^{*}} f\nnum\\&&+\int^{1}_{0}(\Phi^{-1})^{(1,s)^{*}}_{\epsilon Z(t,z(s))}\big(\hat{f}-f\big)ds\Big]\circ z(t)\nnum\\&\doteq&\epsilon H(\epsilon,t,z(t)). \EN

Similar to our analysis in the proof of Theorem \ref{ta}, one can see that $H(0,t,x)=\hat{f}(x)$ where by the construction above, $H$ is smooth with respect to $\epsilon$. By applying the Taylor expansion with remainder we have
\EQ H(\epsilon,t,x)=\hat{f}(x)+\epsilon h(x,\zeta,t),\EN 
where $h(x,\zeta,t)=\frac{\partial }{\partial \epsilon}H(\epsilon,t,x)|_{\epsilon=\zeta}$ and $\zeta\in[0,\epsilon]$.
We note that $H(\epsilon,t,x)$ is periodic with respect to time since $f(x,t)$ and $Z(t,x)$ are both T-periodic. Hence, $h(x,\zeta,t)$ is a T-periodic vector field on $M$. 
Periodicity of $Z$ with respect to $\lambda$ and its continuity with respect to $(\lambda,t)$ give  the compactness of $\hat{D}_{\Phi,\epsilon}$ where 
\EQ \label{bb}\hat{D}_{\Phi,\epsilon}\doteq\bigcup_{\tau\in[0,1],\lambda\in[0,\infty)}\Phi_{\epsilon Z(\lambda,z)}(\tau,0,z_{0}), \hspace{.2cm}\EN 
and similar to the proof of Theorem \ref{ta} we can show that 
\EQ d(\Phi_{\epsilon Z(t,z)}(1,0,z),z)=O(\epsilon),\hspace{.2cm}t\in[t_{0},\infty).\EN
Note that we do not need the compactness of $M$ in order to obtain the statement above since for any initial state $z_{0}$, the state trajectory $\Phi_{\epsilon Z(t,z)}(\tau,0,z_{0}),\tau\in[0,1]$ remains in the compact set $\hat{D}_{\Phi,\epsilon}$. The compactness of $\hat{D}_{\Phi,\epsilon}$ is a direct result of the continuity of $\Phi_{\epsilon Z(\lambda,z)}(\tau,0,z_{0}),\tau\in[0,1]$ with respect to $\tau$ and $\lambda$ (see \cite{mar1}).

The metric triangle inequality implies that
\EQ &&d(\Phi_{\epsilon f(x,t)}(t,t_{0},x_{0}), \Phi_{\epsilon \hat{f}(x)}(t,t_{0},x_{0}))\leq\nnum\\&& d(\Phi_{\epsilon f(x,t)}(t,t_{0},x_{0}),\Phi_{\epsilon Z(t,x)}^{(1,0)}\circ \Phi_{\epsilon f(x,t)}(t,t_{0},x_{0}))\nnum\\&&+ d(\Phi_{\epsilon Z(t,x)}^{(1,0)}\circ \Phi_{\epsilon f(x,t)}(t,t_{0},x_{0}),\Phi_{\epsilon \hat{f}(x)}(t,t_{0},x_{0})).\nnum\EN
Hence, to demonstrate that the hypothesis of the theorem holds, we need to show that
\EQ d(\Phi_{\epsilon Z(t,x)}^{(1,0)}\circ \Phi_{\epsilon f(x,t)}(t,t_{0},x_{0}),\Phi_{\epsilon \hat{f}(x)}(t,t_{0},x_{0}))=O(\epsilon),\hspace{.2cm}\nnum\\t\in[t_{0},\infty).\nnum\EN
To this end, we analyze the distance variation of the following dynamics:
\EQ \label{kkk}&&\dot{x}(t)=\epsilon \hat{f}(x(t)),\hspace{2.2cm}x(t_{0})=x_{0},\nnum\\&&\dot{z}(t)=\epsilon \hat{f}(z(t))+\epsilon^{2} h(z,\zeta,t),\hspace{.2cm}z(t_{0})=x_{0}.\EN

Rescaling time via $\omega=\epsilon t$ in (\ref{kkk}) yields
\EQ &&\frac{dx}{d\omega}=\hat{f}\left(x(\frac{\omega}{\epsilon})\right),\hspace{2.3cm}x(\omega_{0})=x_{0},\nnum\\&&\frac{dz}{d\omega}= \hat{f}\left(z(\frac{\omega}{\epsilon})\right)+\epsilon h(z,\zeta,\frac{\omega}{\epsilon}),\hspace{.2cm}z(\omega_{0})=x_{0}.\EN
Without loss of generality we assume positive definiteness and negative definiteness of $v$ and  $\mathfrak{L}_{\hat{f}}v$ are both defined in the same neighborhood $\mathcal{V}_{\bar{x}}\subset \mathcal{U}_{\bar{x}}$ of $\bar{x}$, where $(\mathcal{U}_{\bar{x}},\psi)$ is the local coordinate system around $\bar{x}$. Otherwise we employ the intersection of the corresponding neighborhoods to perform all the analyses above.
Hence, by Lemma \ref{lc}, there exists $\mathcal{N}_{b}(\bar{x})\subset\mathcal{V}_{\bar{x}}$ such that $\mathcal{N}_{b}(\bar{x})$ is compact. Continuity of solutions and negativity  of $\mathfrak{L}_{\hat{f}}v$ together imply that 
\EQ x\in int(\mathcal{N}_{b}(\bar{x}))\Rightarrow \Phi_{\hat{f}}(\omega,\omega_{0},x)\in \mathcal{N}_{b}(\bar{x}),\hspace{.2cm}\omega\in[\omega_{0},\infty).\nnum\EN\\ 
Now we show that the integral flow of the nominal perturbed system stays close (in the sense of metric $d$) to the integral flow of the averaged system. 
By linearity of Lie derivatives, we have
\EQ \mathfrak{L}_{\hat{f}(z)+\epsilon h(z,\zeta,\omega)}v&=&\mathfrak{L}_{\hat{f}(z)}v+\epsilon \mathfrak{L}_{h(z,\zeta,\omega)}v\nnum\\&=&dv(\hat{f}(z))+\epsilon dv(h(z,\zeta,\omega)).\EN
Since $\mathfrak{L}_{.}(v(x)):T_{x}M\rightarrow \mathds{R}$ is a bounded linear map, we introduce $||\mathfrak{L}_{.}(v)||$ as the operator norm of $dv$. Then we define
\EQ ||\mathfrak{L}_{.}(v)||_{\Omega}\doteq\sup_{x\in \Omega}||\mathfrak{L}_{.}(v(x))||.\EN

Since $\mathcal{N}_{b}(\bar{x})$ contains a neighborhood of $\bar{x}$, applying the Shrinking Lemma  (see \cite{Lee4}) implies the existence of a prempact neighborhood $\mathcal{W}_{\bar{x}}$ such that
\EQ \mathcal{W}_{\bar{x}}\subset \overline{\mathcal{W}_{\bar{x}}}\subset \mathcal{N}_{b}(\bar{x}).\EN 
Then, $M-\mathcal{W}_{\bar{x}}$ is a closed set and $\mathcal{N}_{b}(\bar{x})\bigcap\left( M-\mathcal{W}_{\bar{x}}\right)\subset \mathcal{N}_{b}(\bar{x})$ is a compact set, where 
\EQ \mathfrak{L}_{\hat{f}}v(z)|_{z\in \mathcal{N}_{b}(\bar{x})\bigcap (M-\mathcal{W}_{\bar{x}})}<0.\EN 

 Define  
\EQ \mathfrak{M}\doteq\sup_{z\in \mathcal{N}_{b}(\bar{x})\bigcap (M-\mathcal{W}_{\bar{x}})}\mathfrak{L}_{\hat{f}}v(z)<0.\EN
Consequently 
 
\EQ &&\mathfrak{L}_{\hat{f}(z)}v+\epsilon \mathfrak{L}_{h(z,\zeta,\omega)}v\leq \mathfrak{M}+\epsilon||\mathfrak{L}_{.}v||_{\mathcal{N}_{b}(\bar{x})\bigcap (M-\mathcal{W}_{\bar{x}})}\times\nnum\\&&  \hspace{4cm}||h(z,\zeta,\omega)||_{g},\nnum\\&&z\in\mathcal{N}_{b}(\bar{x})\bigcap (M-\mathcal{W}_{\bar{x}}), \epsilon\in(0,\epsilon_{1}],\omega\in [\omega_{0},\infty),\EN
Since $h$ is periodic with respect to $\omega$ and smooth with respect to $z,\zeta$, there exists $\epsilon_{0}>0$ such that
\EQ  \mathfrak{M}+\epsilon||\mathfrak{L}_{.}v||_{\mathcal{N}_{b}(\bar{x})\bigcap( M-\mathcal{W}_{\bar{x}})} ||h(z,\zeta,\omega)||_{g}\leq 0,\hspace{.2cm}\nnum\\z\in \mathcal{N}_{b}(\bar{x})\bigcap (M-\mathcal{W}_{\bar{x}}), \epsilon\in(0,\epsilon_{0}],\omega\in [\omega_{0},\infty).\EN
Note that $z\in int(\mathcal{N}_{b}(\bar{x}))$ implies that either $z\in \mathcal{N}_{b}(\bar{x})\bigcap (M-\mathcal{W}_{\bar{x}})$ or $z\in \mathcal{W}_{\bar{x}}$. In the first case, $\mathfrak{L}_{\hat{f}(z)}v+\epsilon \mathfrak{L}_{h(z,\zeta,\omega)}v\leq 0$, so that $\Phi_{\hat{f}+\epsilon h}(\omega,\omega_{0},z)\in \mathcal{N}_{b}(\bar{x})$. In the second case $\Phi_{\hat{f}+\epsilon h}(\omega,\omega_{0},z)$ either stays in $\mathcal{W}_{\bar{x}}$ or enters $\mathcal{N}_{b}(\bar{x})\bigcap M-\mathcal{W}_{\bar{x}}$ and consequently it stays in $\mathcal{N}_{b}(\bar{x})$.

Hence, for the interval of existence of solutions $[\omega_{0},\omega_{f})$, we have $\Phi_{\hat{f}+\epsilon h}(\omega,\omega_{0},z)\in \mathcal{N}_{b}(\bar{x}) ,\hspace{.2cm} z\in int(\mathcal{N}_{b}(\bar{x}))$. Applying the Escape Lemma (see \cite{Lee2}) gives  $\omega_{f}=\infty$. That is it has been shown that $\Phi_{\hat{f}+\epsilon h}(\omega,\omega_{0},z)\in \mathcal{N}_{b}(\bar{x}) ,\hspace{.2cm} z\in int(\mathcal{N}_{b}(\bar{x}))$  is bounded in the sense of being trapped in the compact set $\mathcal{N}_{b}(\bar{x})$. 

If the initial state $x_{0}\in \mathcal{N}_{b}(\bar{x})$, then by the statement above $\Phi_{\hat{f}+\epsilon h}(\omega,\omega_{0},x_{0})\in \mathcal{U}_{\bar{x}},\hspace{.2cm}\omega\in[\omega_{0},\infty)$, where $(\mathcal{U}_{\bar{x}},\psi)$ is the local coordinate chart around $\bar{x}$ and (with no loss of generality) we assume $\psi(\bar{x})=0\in \mathds{R}^{n}$.

The uniform normal neighborhood of $\bar{x}\in M$ with respect to $\mathcal{U}_{\bar{x}}$ is denoted by $\mathcal{U}^{n}_{\bar{x}}$ (its existence is guaranteed by Lemma \ref{un}). Consider a  geodesic ball of radius $\delta$ where $\exp_{\bar{x}}(B_{\delta}(0))\subset \mathcal{U}^{n}_{\bar{x}}$.   By definition, $\exp_{\bar{x}}(B_{\delta}(0))$ is an open set containing $\bar{x}$ in the topology of $M$. Therefore by Lemma \ref{lc} one can shrink $b$ to $\acute{b}, 0<\acute{b}\leq b,$  such that  $\mathcal{N}_{\acute{b}}(\bar{x})\subset \exp_{\bar{x}}(B_{\delta}(0))$ ($v$ is locally positive and smooth). Employing the results of \cite{Pet}, Section 5.6, we know that  the distance function $d(\cdot,\bar{x})$ is given locally in $\mathcal{U}^{n}_{\bar{x}}$ by
\EQ d(x,\bar{x})=\left(\sum^{n}_{i=1}x^{2}_{i}\right)^{\frac{1}{2}},\EN
which is the Euclidean distance function and hence in the normal coordinate system the convergence in the topology of $M$ will be same as the convergence in the Euclidean topology. The vector space $T_{x}M$ is a finite dimensional normed vector space therefore $\overline{B}_{\delta}(0)\subset T_{\bar{x}}M$ is compact and consequently $\exp_{\bar{x}}(\overline{B}_{\delta}(0))\subset M$ is a compact set ($\exp$ is a local diffeomorphism). 
Let us replace the Riemannian metric $g$ with the standard Euclidean metric on $\mathcal{N}_{\acute{b}}$. Smoothness of $\hat{f}$ and compactness of $\exp_{\bar{x}}(\overline{B}_{\delta}(0))$ together imply that the Jacobian matrix $\frac{\partial\hat{f}}{\partial x}$ is bounded on $\mathcal{N}_{b}$ and hence the conditions of the Converse Lyapunov Theorem (see \cite{Kha}, Theorem 4.14) are satisfied. 
Since $\bar{x}$ is exponentially stable, invoking the results of \cite{Kha} (Theorems 4.14 and 9.1) implies that there exists a parameter $0<\beta$ which is independent of $||h(z,\zeta,\omega)||_{e}$ such that
\EQ \label{luc}||x(\omega)-z(\omega)||_{e}\leq \epsilon\beta K,\EN 
where  $||.||_{e}$ is the Euclidean norm of $\mathds{R}^{n}$, and 
\EQ K&\doteq&\sup_{z\in \mathcal{N}_{b}(\bar{x}), \omega\in[\omega_{0},\infty)}||h(z,\zeta,\omega)||_{e}\nnum\\&=&\sup_{z\in \mathcal{N}_{b}(\bar{x}), \omega\in[\omega_{0},\omega_{0}+ T]}||h(z,\zeta,\omega)||_{e}.\EN  
 
 The vector space $T_{x}M$ is scalable with respect to $||.||_{e}$ and $||.||_{g}$, i.e. there exist $0<\lambda_{1}\leq \lambda_{2}$ such that $\lambda_{1}||X||_{e}\leq ||X||_{g}\leq \lambda_{2}||X||_{e},\hspace{.2cm} X\in T_{x}M$.
  Continuity of $v$ implies that $\mathcal{N}_{\acute{b}}(\bar{x})$ is closed and compact, with the latter following as $\mathcal{N}_{\acute{b}}(\bar{x})\subset \mathcal{N}_{b}(\bar{x})$ and $\mathcal{N}_{b}(\bar{x})$ is compact. Hence, without loss of generality, we  assume that $\mathcal{N}_{b}(\bar{x})\subset \exp_{\bar{x}}(B_{\delta}(0))$ and consequently $\Phi_{\hat{f}+\epsilon h}(\omega,\omega_{0},z)\in \exp_{\bar{x}}(B_{\delta}(0)),\hspace{.2cm}\omega\in[\omega_{0},\infty), z\in \mathcal{N}_{b}(\bar{x}) $. By scaling  the Euclidean and Riemannian metrics inside $\exp_{\bar{x}}(\overline{B}_{\delta}(0))$ we have (for the scaling procedure see \cite{Lee3}, Lemma 5.12)
\EQ\lambda_{1}||X||_{e}\leq ||X||_{g}\leq \lambda_{2}||X||_{e},\hspace{.2cm}X\in T_{x}M, x\in  \exp_{\bar{x}}(\overline{B}_{\delta}(0)).\nnum\\ \EN
Now in the Euclidean metric $||.||_{e}$ consider a  smooth straight line parametrized by time  $\gamma_{1}:[0,1]\rightarrow \mathds{R}^{n}$ such that $\gamma_{1}(0)=x(\omega)$ and $\gamma_{1}(1)=z(\omega)$. 

The results of \cite{Pet}, Corollary 5.3, ensure that when $\exp$ is  diffeomorphic on its image then  the Euclidean distance ball and the geodesic ball are identical sets on $M$. Therefore employing (\ref{luc}) implies that choosing $\hat{\epsilon}<min\{\epsilon_{0},\epsilon_{1}\}$ small enough guarantees the closeness of  $x(\omega)$ and $z(\omega)$  in the sense that $\gamma_{1}\subset \exp_{\bar{x}}(B_{\delta}(0))$. Hence,
\EQ \label{tay}&&d(x(\omega),z(\omega))\leq \ell(\gamma_{1})= \int^{1}_{0}||\dot{\gamma_{1}}(\tau)||_{g}d\tau\leq\nnum\\&& \lambda_{2}\int^{1}_{0}||\dot{\gamma_{1}}(\tau)||_{e}d\tau=\lambda_{2}||x(\omega)-z(\omega)||_{e}\leq \epsilon \lambda_{2} \beta K,\nnum\\&&\omega\in [\omega_{0},\infty), \EN
which completes the proof by choosing $\mathcal{N}_{\bar{x}}=int(\mathcal{N}_{b}(\bar{x}))$.\\\hspace*{8cm}\qed
\end{pf}
The exponential stability requirement of Theorem \ref{taa} can be relaxed to local asymptotic stability. This yields a closeness of solutions result of a similar form to Theorem \ref{taa}, but with a weaker
implied property. A version of this result is obtained for dynamical systems with external disturbances in \cite{nes}, Theorem 1. In \cite{Peut}, for a special case of homogeneous dynamical systems and under some technical hypotheses, it has been shown that the asymptotic stability of the average system implies the asymptotic stability of the nominal system.    
\begin{theorem}
\label{taaa}
For a smooth $n$ dimensional Riemannian manifold $M$, let $f\in \mathfrak{X}(M\times \mathds{R})$ be a $T$-periodic smooth vector field and assume  the nominal and averaged vector fields are both complete for $\epsilon\in(0,\epsilon_{1}],\hspace{.2cm}0< \epsilon_{1}$. Suppose the averaged dynamical system has a locally asymptotically stable equilibrium $\bar{x}\in M$ such that there exists a Lyapunov function $v:M\rightarrow \mathds{R}$  with $\mathfrak{L}_{\hat{f}}v$ is locally negative-definite around $\bar{x}$. Then, for every $0<\delta$, there exists a neighborhood $\mathcal{N}_{\bar{x}}$ and $\hat{\epsilon}\leq\epsilon_{1}$, such that   
\EQ d(\Phi_{\epsilon f}(t,t_{0},x_{0}),\Phi_{\epsilon \hat{f}}(t,t_{0},x_{0}))\leq \delta,\hspace{.2cm}\nnum\\\epsilon\in(0,\hat{\epsilon}],x_{0}\in \mathcal{N}_{\bar{x}},\hspace{.2cm} t\in[t_{0},\infty).\EN
\end{theorem}
\begin{pf}
Without loss of generality  assume $\delta$ is  small enough so that $\exp$ is a diffeomorphism on $\exp_{\bar{x}}(B_{\frac{\delta}{4}}(0))$ and $\exp_{\bar{x}}(B_{\frac{\delta}{4}}(0))\subset \mathcal{U}^{n}_{\bar{x}}$ where $B_{\frac{\delta}{4}}(0)\in T_{\bar{x}}M$ ($\mathcal{U}^{n}_{\bar{x}}$ is the uniform normal neighborhood around $\bar{x}$). Following the steps of the proof of Theorem \ref{taa} it can be shown that $\Phi_{\hat{f}+\epsilon h}(\omega,\omega_{o},z)\in\mathcal{N}_{b}(\bar{x})  ,\hspace{.2cm} z\in int(\mathcal{N}_{b}(\bar{x}))$ where $\mathcal{N}_{b}(\bar{x})$ is a connected compact sublevel set of the Lyapunov function $v$ such that $\mathcal{N}_{b}(\bar{x})\subset \exp_{\bar{x}}(B_{\frac{\delta}{4}}(0))$. Now let us consider the Euclidean metric instead of the Riemannian one on $\exp_{\bar{x}}(B_{\frac{\delta}{4}}(0))$. Since $\exp_{\bar{x}}(B_{\frac{\delta}{4}}(0))\subset \mathcal{U}^{n}_{\bar{x}}$, employing the results of \cite{Pet}, Corollary 5.3, \cite{Lee3}, Proposition 5.11, implies that the geodesic balls and Euclidean balls are identical, while smoothness of $\hat{f}$ implies the boundedness of $\frac{\partial \hat{f}}{\partial x}$ on the compact set $\exp_{\bar{x}}(\overline{B}_{\frac{\delta}{4}}(0))$. Combining the results of \cite{Kha}, Theorem 4.16 and Lemma 9.3 together implies that  there exists $\epsilon_{0}>0$ such that
\EQ \label{hel}d(\Phi_{\hat{f}+\epsilon h}(\omega,\omega_{0},x),\bar{x})\leq \rho(\epsilon K),\hspace{.2cm}  \omega_{0}\leq \omega<\infty,\epsilon\in(0,\epsilon_{0}],\nnum\\\EN
where $K\doteq \sup_{x\in \mathcal{N}_{b}(\bar{x}),\omega\in [\omega_{0},\omega_{0}+T]}||h(x,\zeta,\omega)||_{e},$
and $\rho$ is a strictly increasing continuous function satisfying $\rho(0)=0$. Note that there exists a class $\mathcal{K}\mathcal{L}$ function (for the definition of $\mathcal{K}\mathcal{L}$ functions see \cite{Kha}, Section 4.4) in the statement of Lemma 9.3 in \cite{Kha}  which bounds the state trajectory up to a specified time $t$. Since this $\mathcal{K}\mathcal{L}$ function is decreasing with respect to time and its construction only depends on the average system then we can choose $\epsilon_{0}$ sufficiently small such that (\ref{hel}) holds for all $\omega\geq \omega_{0}$. 

Now by the continuity of $\rho$, we can choose $\epsilon_{0}$ sufficiently small such that $d(\Phi_{\hat{f}+\epsilon h}(\omega,\omega_{0},x),\bar{x})\leq \rho(\epsilon K)\leq \frac{\delta}{4}, \forall\epsilon\in(0,\epsilon_{0}]$. Selecting the initial condition $x_{0}\in \mathcal{N}_{b}(\bar{x})\subset \exp_{\bar{x}}(B_{\frac{\delta}{4}}(0)) $ guarantees that the state trajectory of the average system does not exit $\mathcal{N}_{b}(\bar{x})$. Hence, 
\EQ &&d(\Phi_{\hat{f}+\epsilon h}(\omega,\omega_{0},x_{0}),\Phi_{\hat{f}}(\omega,\omega_{0},x_{0}))\leq\nnum\\&& d(\Phi_{\hat{f}+\epsilon h}(\omega,\omega_{0},x_{0}),\bar{x})+d(\Phi_{\hat{f}}(\omega,\omega_{0},x_{0}),\bar{x})\leq \frac{\delta}{2},\hspace{.2cm}\nnum\\&&\forall x_{0}\in \mathcal{N}_{b}(\bar{x}).\EN
Therefore, we choose $\hat{\epsilon}\leq min\{\epsilon_{0},\epsilon_{1}\}$, so that 
\EQ &&d(\Phi_{\epsilon f(x,t)}(t,t_{0},x_{0}), \Phi_{\epsilon \hat{f}(x)}(t,t_{0},x_{0}))\leq\nnum\\&& d(\Phi_{\epsilon f(x,t)}(t,t_{0},x_{0}),\Phi_{\epsilon Z(t,x)}^{(1,0)}\circ \Phi_{\epsilon f(x,t)}(t,t_{0},x_{0}))+ \nnum\\&&d(\Phi_{\epsilon Z(t,x)}^{(1,0)}\circ \Phi_{\epsilon f(x,t)}(t,t_{0},x_{0}),\Phi_{\epsilon \hat{f}(x)}(t,t_{0},x_{0}))\leq\nnum\\&& O(\epsilon)+\frac{\delta}{2}\leq \delta,\hspace{.2cm}\epsilon\in(0,\hat{\epsilon}].\nnum\EN
Hence, the statement of the theorem follows for $\mathcal{N}_{\bar{x}}=int(\mathcal{N}_{b}(\bar{x}))$ and $\hat{\epsilon}$.
\qed\end{pf}
\section{Almost global stability and infinite horizon averaging on compact Riemannian manifolds}
Now we focus on the analysis of the closeness of solutions for dynamical systems evolving on compact Riemannian manifolds where the average system is almost globally stable. The notion of almost global stability is defined below. We note that, due to the non-contractibility of compact manifolds, there exists no smooth vector field which globally asymptotically stabilizes an equilibrium on a compact configuration manifold, see \cite{San, Mai}. 
  \begin{definition}[\cite{San, Mai}]
  \label{nic}
For the dynamical system $\dot{x}(t)=f(x), \hspace{.2cm}f:M\rightarrow TM$, an equilibrium $\bar{x}\in M$ is almost globally asymptotically/exponentially stable if there exists an open $\mathcal{U}_{\bar{x}}$ dense in $M$ such that for all $t_{0}\in \mathds{R}$\\
 
(i):(almost globally asymptotically) $\bar{x}$  is Lyapunov stable on $M$ and 
\EQ \forall x(t_{0})\in \mathcal{U}_{\bar{x}},\hspace{.2cm}\lim_{t\rightarrow \infty}\Phi_{f}(t,t_{0},x(t_{0}))=\bar{x}. \EN

(ii):(almost globally exponentially) if $\bar{x}$ is almost globally asymptotically stable and \EQ \label{tas}d(\Phi_{f}(t,t_{0},x(t_{0})),\bar{x})\leq &&kd(x(t_{0}),\bar{x})\exp(-\lambda(t-t_{0})),\nnum\\&& k,\lambda\in \mathds{R}_{>0},x(t_{0})\in \mathcal{U}_{\bar{x}}.\EN \qed
\end{definition}
 The following Theorems specify closeness of solutions on an infinite time horizon for systems evolving on compact Riemannian
manifolds.
\begin{theorem}
\label{ttaa}
For a smooth $n$ dimensional compact Riemannian manifold $M$, let $f\in \mathfrak{X}(M\times\mathds{R})$ be a $T$-periodic smooth vector field. Suppose $\bar{x}\in M$ is almost globally  exponentially stable on $M$ for the average dynamical system $\hat{f}$ and there exists a Lyapunov function $v:M\rightarrow \mathds{R}$  such that $\mathfrak{L}_{\hat{f}}v$ is locally negative-definite around $\bar{x}$. Then, there exist a dense open set $\mathcal{U}_{\bar{x}}\subset M$ and $\hat{\epsilon}>0$ such that   
\EQ d(\Phi_{\epsilon f}(t,t_{0},x_{0}),\Phi_{\epsilon \hat{f}}(t,t_{0},x_{0}))=O(\epsilon),\hspace{.2cm}\nnum\\\epsilon\in(0,\hat{\epsilon}],\forall x_{0}\in \mathcal{U}_{\bar{x}},\hspace{.2cm} t\in[t_{0},\infty).\EN
\end{theorem}
\begin{pf}
The proof follows via Lemma \ref{lp} and Theorem \ref{taa}.
First we note that, since $\hat{f}$ is almost globally exponentially stable, by Definition \ref{nic} there exists $\mathcal{U}_{\bar{x}}$ such that (\ref{tas}) holds. Since $\mathcal{U}_{\bar{x}}$ is open in the topology of $M$ then $M-\mathcal{U}_{\bar{x}}$ is closed and closed subsets of compact sets are all compact, see \cite{Lee4}. Hence, there exists  a neighborhood $ U^{1}_{\bar{x}}\subset M$ such that $U^{1}_{\bar{x}}\bigcap (M-\mathcal{U}_{\bar{x}})=\emptyset$. Otherwise, $\bar{x}\in M-\mathcal{U}_{\bar{x}}$, or $\bar{x}$ is a limit point of $M-\mathcal{U}_{\bar{x}}$. Since $M-\mathcal{U}_{\bar{x}}$ is closed, it follows that $\bar{x}\in M-\mathcal{U}_{\bar{x}}$ which contradicts the fact that $\hat{f}$ is almost globally exponentially stable on $\mathcal{U}_{\bar{x}}$. In the time scaled variable $\omega=\epsilon t$, the exponential stability of $\hat{f}$ implies that there exist $\tau>0$ such that $\Phi_{\hat{f}}(\omega,\omega_{0},x_{0})\in U^{1}_{\bar{x}},\hspace{.2cm}\omega\geq\tau$, see Figure \ref{mm}, 
\begin{figure}
\begin{center}
\hspace*{-.25cm}\includegraphics[scale=.3]{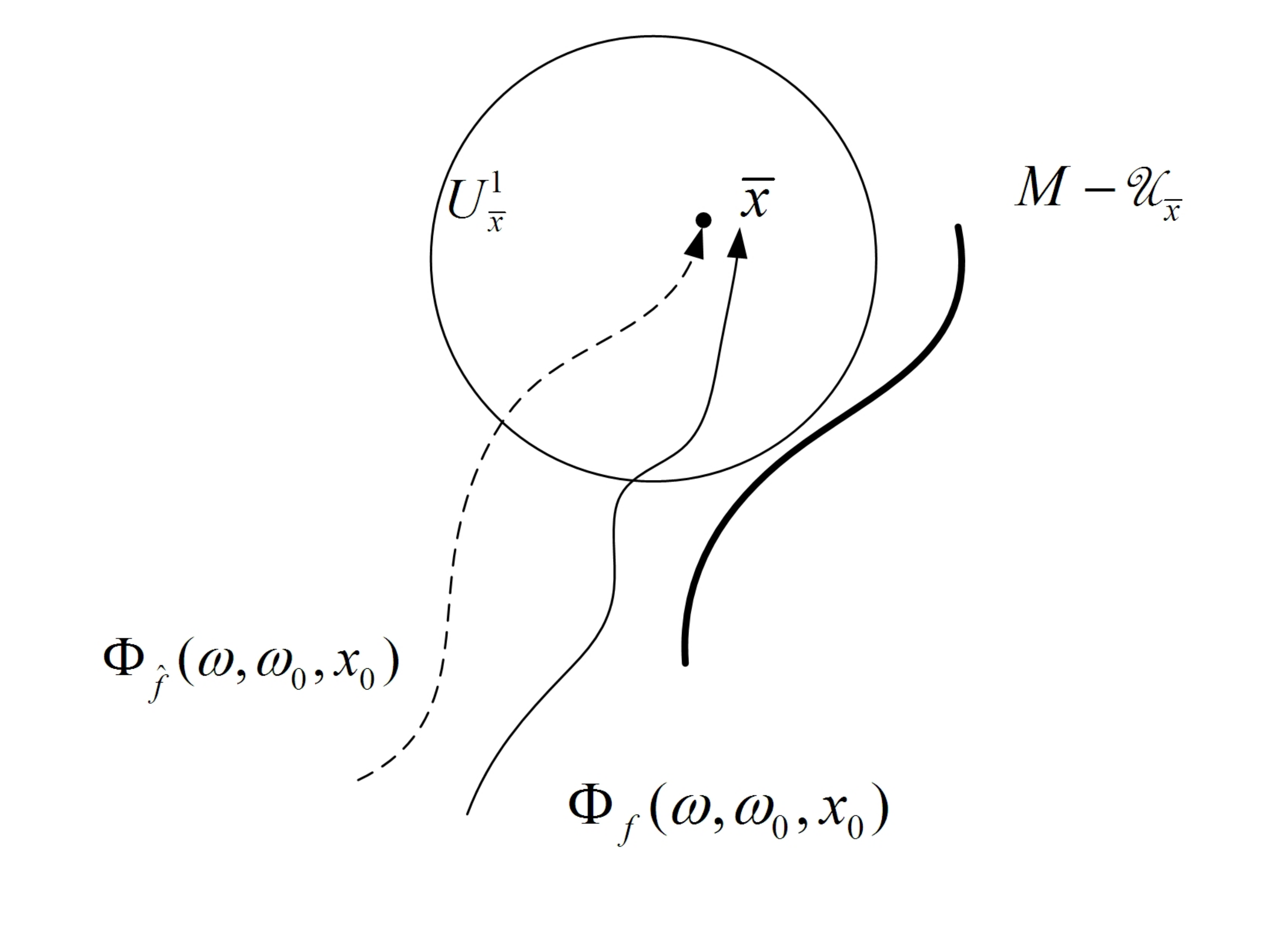}
 \caption{State trajectories of the nominal and averaged systems.} 
     \label{mm}
      \end{center}
   \end{figure}  
and also continuity of $\Phi_{\hat{f}}(.,\omega_{0},x_{0})$ gives the compactness of $\bigcup_{\omega\in[\omega_{0},\tau]}\Phi_{\hat{f}}(\omega,\omega_{0},x_{0})$ in $M$. The distance function $d$ on $M$ is continuous with respect to both of its arguments, so that by Lemma \ref{lp}, on compact time intervals (compactness in $\omega$ gives $t-t_{0}=O(\frac{1}{\epsilon})$) we can select $\epsilon_{1}$ and $U^{1}_{\bar{x}}$ as small as 
\EQ \left(\bigcup_{\omega\in[\omega_{0},\tau]}\Phi_{ f}(\omega,\omega_{0},x_{0})\right)\bigcap (M-\mathcal{U}_{\bar{x}})=\emptyset.\EN
 As presented in the proof of Theorem \ref{taa}, we have the following time rescaled equations:
\EQ &&\frac{dx}{d\omega}=\hat{f}\left(x(\frac{\omega}{\epsilon})\right),\hspace{2.3cm}x(\omega_{0})=x_{0},\nnum\\&&\frac{dz}{d\omega}= \hat{f}\left(z(\frac{\omega}{\epsilon})\right)+\epsilon h(z,\zeta,\frac{\omega}{\epsilon}),\hspace{.2cm}z(\omega_{0})=x_{0}.\EN
Employing the results of  Lemma \ref{lp}, on a compact interval of time $\omega_{1}-\omega_{0}$, we can shrink $\epsilon_{1}$  as  
\EQ \label{jus}d\left(\Phi_{\hat{f}}(\omega,\omega_{0},x_{0}),\Phi_{\hat{f}+\epsilon h}(\omega,\omega_{0}, x_{0})\right)=O(\epsilon),\hspace{.2cm}\nnum\\\omega\in[\omega_{0},\omega_{1}],\epsilon\in(0,\epsilon_{1}].\EN 
Assume for $\delta >0$ that $\exp_{\bar{x}}(B_{\delta}(0))\subset\mathcal{U}^{n}_{\bar{x}}$, where $\mathcal{U}^{n}_{\bar{x}}$ is a uniform normal neighborhood of $\bar{x}$ and choose $\omega_{1}=\tau$ and $U^{1}_{\bar{x}}=int(\mathcal{N}_{b}(\bar{x}))\subset \exp_{\bar{x}}(B_{\delta}(0))$ where $\mathcal{N}_{b}(\bar{x})$ is a compact connected sublevel set of $v$. We note that the existence of $\omega_{1}$ does not guarantee the existence of an entry time $t_{1}$ for the unscaled dynamical system since the smaller we choose $\epsilon$, the larger time it takes for the state trajectory to enter $\exp_{\bar{x}}(B_{\delta}(0))$. Now we show that $d\big(\Phi_{\hat{f}}(\omega,\omega_{1},x(\omega_{1})),\Phi_{\hat{f}+\epsilon h}(\omega,\omega_{1}, z(\omega_{1}))\big)=O(\epsilon),\hspace{.2cm}\omega_{1}<\omega$. 

Following the proof of Theorem \ref{taa} and employing the results of \cite{Kha}, and Theorems 4.14 and 9.1, we have
\EQ ||x(\omega)-z(\omega)||_{e}&\leq& k\exp(-\gamma(\omega-\omega_{1}))||x(\omega_{1})-z(\omega_{1})||_{e}\nnum\\
&&+\beta \epsilon K,\omega\in[\omega_{1},\omega_{2}),\EN
for some parameters $k,\beta,\gamma$ which are independent of $\epsilon$ and $K$ ($K$ is defined as per the proof of Theorem \ref{taa}) and $z(\omega)\in\mathcal{U}^{n}_{\bar{x}},\hspace{.2cm}\omega\in[\omega_{1},\omega_{2})$. It remains to show that the Euclidean distance can be scaled by the Riemannian distance. In the last part of the proof of Theorem \ref{taa} we have shown that the Riemannian distance can be bounded above by the Euclidean distance. Similar to the scaling procedure presented in the proof of Theorem \ref{taa} we can show there exist $\hat{\lambda}_{1},\hat{\lambda}_{2}\in \mathds{R}_{>0}$ such that
\EQ\hat{\lambda}_{1}||X||_{g}\leq ||X||_{e}\leq \hat{\lambda}_{2}||X||_{g},\hspace{.2cm}X\in T_{x}M, x\in  \mathcal{N}_{b}(\bar{x}),\nnum\\ \EN
where $\mathcal{N}_{b}(\bar{x}) $ is a compact set.
By Lemma \ref{lp}, we select $\hat{\epsilon}\leq\epsilon_{1}$ sufficiently small such that $x(\omega_{1}),z(\omega_{1})\in \mathcal{U}_{x(\omega_{1})}\subset\mathcal{N}_{b}(\bar{x}),\hspace{.2cm}\forall\epsilon\in(0,\hat{\epsilon}]$, where $\mathcal{U}_{x(\omega_{1})}=\{x\in int(\mathcal{N}_{b}(\bar{x})), ||x-x(\omega_{1})||_{e}<\rho\}$ is an open set for a sufficiently small $\rho\in \mathds{R}_{>0}$. Now consider an arbitrary piecewise smooth curve $\gamma:[0,1]\rightarrow M$ connecting $x(\omega_{1})$ and $z(\omega_{1})$. Suppose that $\gamma\subset\mathcal{U}_{x(\omega_{1})}\subset\mathcal{N}_{b}(\bar{x}) $ then
\EQ ||x(\omega_{1})-z(\omega_{1})||_{e}\leq \int^{1}_{0}||\dot{\gamma}(\tau)||_{e}d\tau&\leq&\hat{\lambda}_{2}\int^{1}_{0}||\dot{\gamma}(\tau)||_{g}d\tau.\nnum\EN
In the case $\gamma\nsubset\mathcal{U}_{x(\omega_{1})},$ there exists a hitting time $\tau_{1}$ such that $\gamma_{[0,\tau_{1})}\subset \mathcal{U}_{x(\omega_{1})}$ and $\gamma(\tau_{1})\in \overline{\mathcal{U}}_{x(\omega_{1})}$. Since $x(\omega_{1}),z(\omega_{1})\in \mathcal{U}_{x(\omega_{1})}$ and  
\EQ||x(\omega_{1})-z(\omega_{1})||_{e}&\leq&\rho\leq\int^{\tau_{1}}_{0}||\dot{\gamma}(\tau)||_{e}d\tau\nnum\\&\leq&\hat{\lambda}_{2}\int^{\tau_{1}}_{0}||\dot{\gamma}(\tau)||_{g}d\tau\leq \hat{\lambda}_{2}\int^{1}_{0}||\dot{\gamma}(\tau)||_{g}d\tau.\nnum\\\EN
Hence, in general, for all piecewise smooth $\gamma$, $||x(\omega_{1})-z(\omega_{1})||_{e}\leq \hat{\lambda}_{2}\int^{1}_{0}||\dot{\gamma}(\tau)||_{g}d\tau$.Taking the infimum of the right hand side of the equation above implies $||x(\omega_{1})-z(\omega_{1})||_{e}\leq \hat{\lambda}_{2}d(\Phi_{\hat{f}}(\omega_{1},\omega_{0},x_{0}),\Phi_{\hat{f}+\epsilon h}(\omega_{1},\omega_{0}, x_{0})).$
Therefore, we can extend $\omega_{2}$ to $\infty$ and
\EQ\label{kay}||x(\omega)-z(\omega)||_{e}=O(\epsilon), \omega\in[\omega_{1},\infty).\EN
The theorem statement follows by (\ref{jus}) and applying the last part of the proof of Theorem \ref{taa}, i.e. (\ref{tay}) to (\ref{kay}).
\qed\end{pf}
The following Theorem specifies closeness of solutions on an infinite time horizon for systems evolving on compact
Riemannian manifolds in the case where the average system is almost globally asymptotically stable.
\begin{theorem}
\label{ttaaa}
For a smooth $n$ dimensional compact Riemannian manifold $M$, let $f\in \mathfrak{X}(M\times\mathds{R})$ be a $T$-periodic smooth vector field. Suppose $\bar{x}\in M$ is almost globally  asymptotically stable on $ M$ for the average dynamical system $\hat{f}$ and there exists a Lyapunov function $v:M\rightarrow \mathds{R}$  with $\mathfrak{L}_{\hat{f}}v$ is locally negative-definite around $\bar{x}$. Then for every $\delta>0,$ there exist $\mathcal{U}_{\bar{x}}\subset M$ and $\hat{\epsilon}>0$, such that   
$d\left(\Phi_{\epsilon f}(t,t_{0},x_{0}),\Phi_{\epsilon \hat{f}}(t,t_{0},x_{0})\right)<\delta ,\hspace{.2cm}\nnum\\\epsilon\in(0,\hat{\epsilon}],\forall x_{0}\in \mathcal{U}_{\bar{x}},\hspace{.2cm} t\in[t_{0},\infty).$
\end{theorem}
\begin{pf}
The proof parallels that of Theorem \ref{ttaa} by employing the results of Theorem \ref{taaa}.
\qed\end{pf}
\subsection{Example 2}
\label{ex2}
Consider the following dynamical system on a torus $\textbf{T}^{2}$. A parametrization of $\textbf{T}^{2}$ is given by
\EQ x(\theta_{1},\theta_{2})&=&(R+r\cos(\theta_{1}))\cos(\theta_{2}),\nnum\\y(\theta_{1},\theta_{2})&=&(R+r\cos(\theta_{1}))\sin(\theta_{2}),\nnum\\z(\theta_{1},\theta_{2})&=&r \sin(\theta_{2}),\hspace{.2cm}\theta_{1},\theta_{2}\in[-\pi,\pi].\EN
The induced Riemannian metric is given by $g_{T^{2}}(\theta_{1}, \theta_{2})\doteq(R+r\cos(\theta_{1}))^{2}d\theta_{2} \otimes d\theta_{2}+r^{2}d\theta_{1} \otimes d\theta_{1},\hspace{.2cm}R=1,r=.5$,
 \begin{figure}
        \begin{center}
\hspace*{-.25cm}\includegraphics[scale=.25]{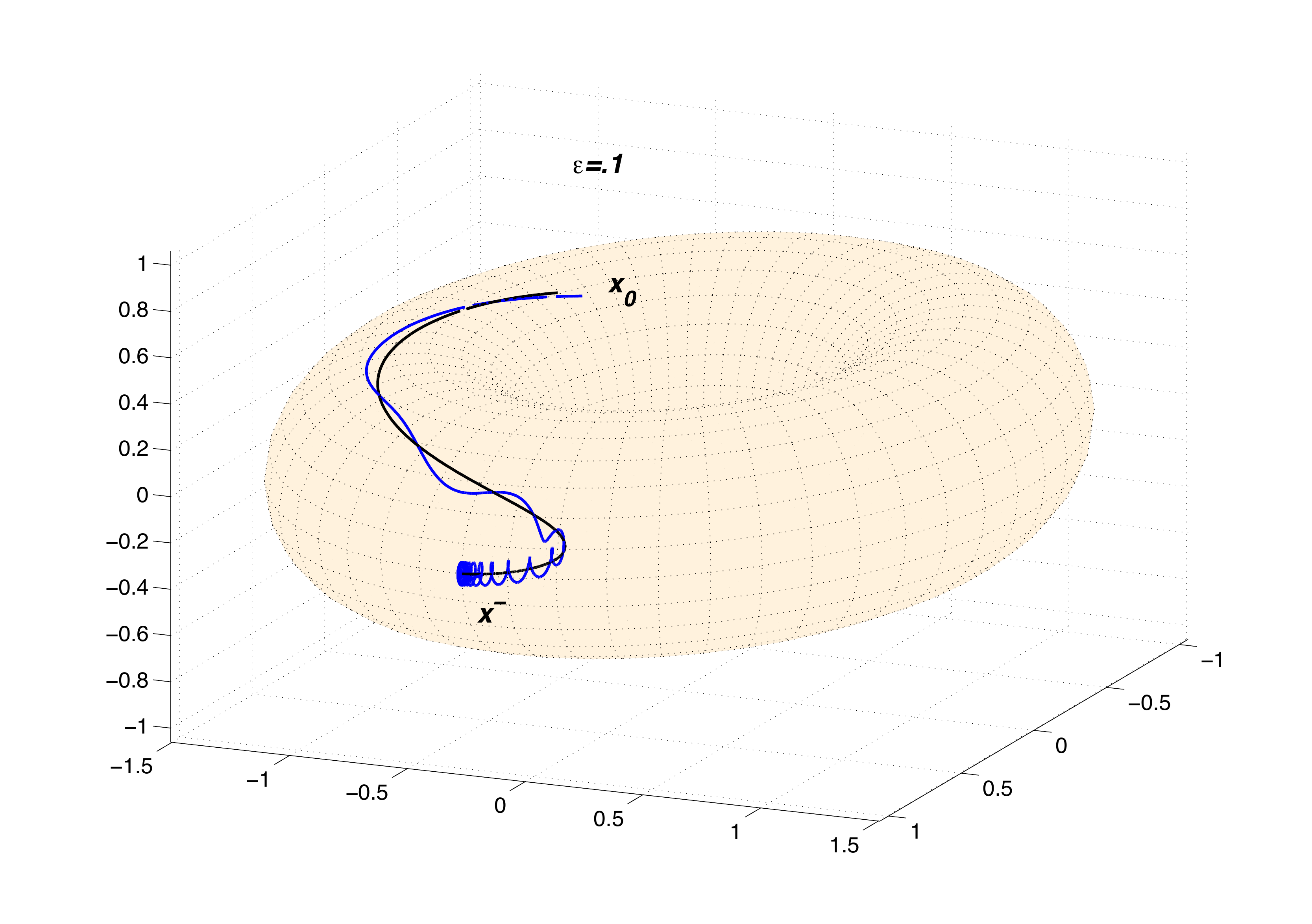}
 \caption{State trajectories on the torus for example 2 (nominal system: blue line, average system: black line), $\epsilon=.1.$ (Section \ref{ex2})}
     \label{f11}
      \end{center}
   \end{figure}  
   \begin{figure}
\begin{center}
\hspace*{-.25cm}\includegraphics[scale=.25]{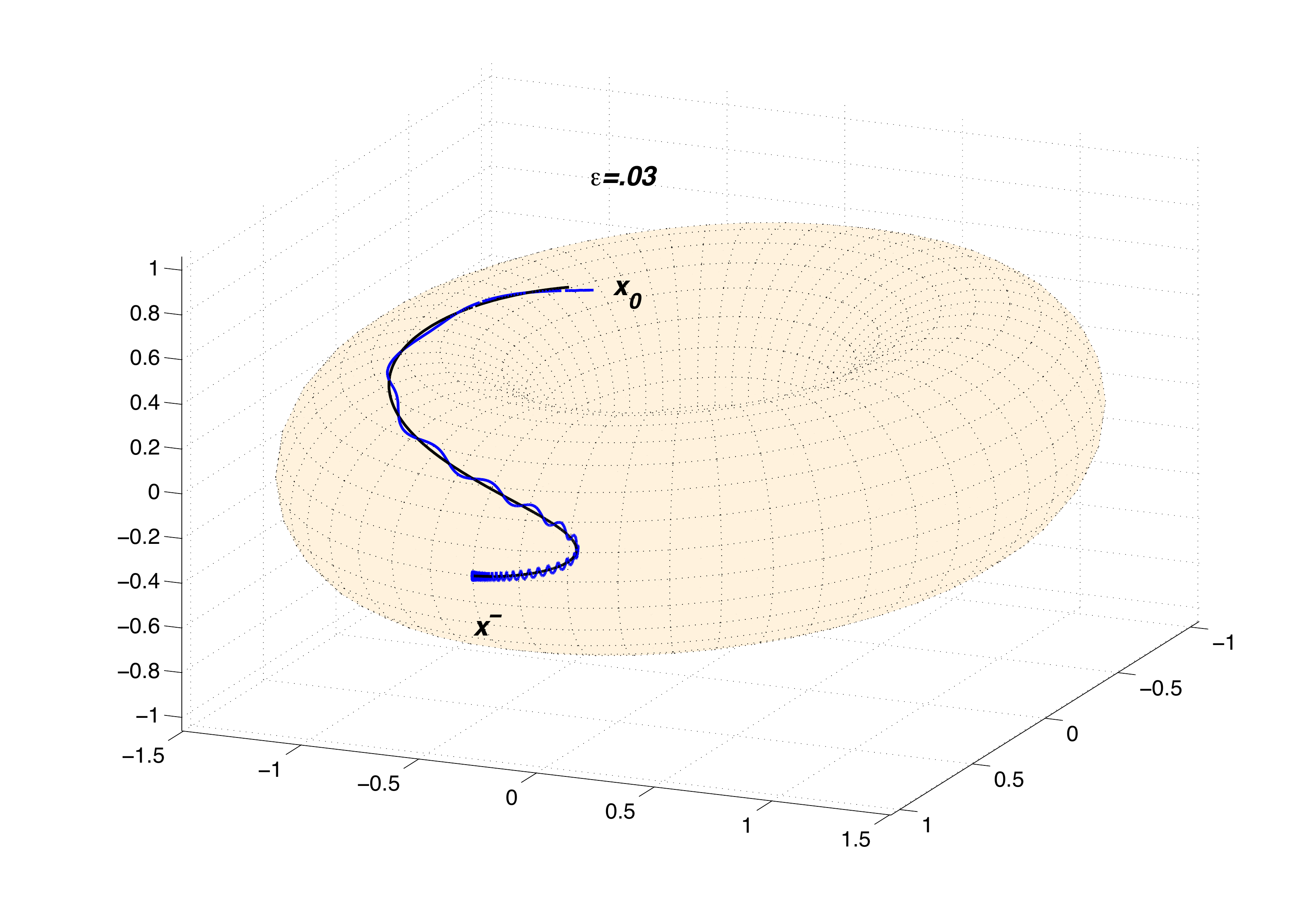}
 \caption{State trajectories on the torus for example 2 (nominal system: blue line, average system: black line), $\epsilon=.03.$ (Section \ref{ex2})}
     \label{f22}
      \end{center}
   \end{figure} 
    
where $\otimes$ is the tensor product, see \cite{Lee2}. 
The dynamical equations are as follows:
\vspace{-.5cm}
\EQ \label{e1}f:\left\{ \begin{array}{ll}
         \dot{\theta_{1}}(t)=\epsilon(-\theta_{1}(t)-\sin(t)),\\
        \dot{\theta_{2}}(t)=\epsilon(\theta_{1}(t)-\theta_{2}(t)).\end{array} \right.\EN
        By applying (\ref{ppp}) to (\ref{e1}), the averaged system is given by
        \EQ\hat{f}:\left\{ \begin{array}{ll}
         \dot{\theta_{1}}(t)=-\epsilon\theta_{1}(t),\\
        \dot{\theta_{2}}(t)=\epsilon(-\theta_{1}(t)-\theta_{2}(t)). \end{array} \right.\EN
   By inspection, the averaged system is locally exponentially stable in a neighborhood of $(0,0)$ for the Euclidean metric on $\textbf{T}^{2}$. By the scaling method of the Riemannian and Euclidean metrics (see \cite{Lee3}), we can show that (\ref{tas}) holds locally around $(0,0)$.      
   Figures \ref{f11} and \ref{f22} show the closeness of solutions for the nominal and averaged systems above for $\epsilon=.1,$ and $.03$ respectively for $t\in[0,\infty)$ as expected by the results of Theorem \ref{taa}. 

\bibliographystyle{plain}        
\bibliography{HSCC1}           


\appendix
\end{document}